\documentclass{amsart}
\usepackage{ifthen,amscd}
\makeatletter
\DeclareRobustCommand{\qbinom}{\genfrac[]\z@{}}
\makeatother
\newcommand{\nk}{\mathbb K}

\newcommand{\spmapright}[1]{\smash{\mathop{
\hbox to 1cm{\rightarrowfill}}\limits^{#1}}}
\newcommand{\lspmapright}[1]{\smash{\mathop{
\hbox to 2cm{\rightarrowfill}}\limits^{#1}}}
\newcommand{\sbmapright}[1]{\smash{\mathop{
\hbox to 1cm{\rightarrowfill}}\limits_{#1}}}

\newcommand{\rmapdown}[1]{\Big\downarrow
\rlap{$\vcenter{\hbox{$\scriptstyle#1$}}$}}

\newtheorem{thm}{Theorem}[subsection]
\newtheorem{lem}[thm]{Lemma}
\newtheorem{prop}[thm]{Proposition}

\newtheorem{cor}[thm]{Corollary}
\theoremstyle{definition}
\newtheorem{defn}[thm]{Definition}

\newtheorem{conj}{Conjecture}[section]
\newtheorem{theorem}[conj]{Theorem}
\newtheorem{ex}[thm]{Example}
\theoremstyle{remark}
\newtheorem{rem}[thm]{Remark}
\numberwithin{equation}{section}
\newenvironment{Ac}%
 {\hspace*{-1.2em}\textbf{Acknowledgment.}\hspace{0.7em}}{}
 {\hspace*{-1.2em}\textbf{Notation.}\hspace{0.7em}}{}
\allowdisplaybreaks
\DeclareMathOperator{\id}{id}
\DeclareMathOperator{\soc}{soc}

\let\top\Top
\DeclareMathOperator{\Hom}{Hom}
\DeclareMathOperator{\Ker}{Ker}
\DeclareMathOperator{\Img}{Im}
\DeclareMathOperator{\Ext}{Ext}
\DeclareMathOperator{\End}{End}
\DeclareMathOperator{\Ob}{Ob}
\let\to\longrightarrow
\let\mapsto\longmapsto
\title[Tensor products of modules over $\overline{U}_q(\mathfrak{sl}_2)$]%
{Indecomposable decomposition of 
tensor products of modules over the restricted quantum universal
enveloping algebra associated to $\boldsymbol{\mathfrak{sl}_2}$}
\author[H.~Kondo and Y.~Saito]{Hiroki KONDO and Yoshihisa Saito}
\address{Graduate School of Mathematical Sciences, 
University of Tokyo, 3-8-1 Komaba,
Meguro-ku, Tokyo 153-8914, Japan.}
\email{HK:donko@ms.u-tokyo.ac.jp\\
YS:yosihisa@ms.u-tokyo.ac.jp}
\begin{document}
\begin{abstract}
We study 
the tensor structure of the category of finite dimensional modules of   
the restricted quantum enveloping algebra 
associated to $\mathfrak{sl}_2$. 
Tensor product decomposition rules for all
indecomposable modules are explicitly given. As a by-product,
it is also shown that the category of finite dimensional modules
of the restricted quantum enveloping algebra 
associated to $\mathfrak{sl}_2$ is {\it not} a braided tensor category.
\end{abstract}
\maketitle
\section{Introduction}
In the representation theory of quantum groups at roots of unity, it is
often assumed that the parameter $q$ is a primitive $n$-th root of unity
where $n$ is a odd prime number. However, there has recently been
increasing interest in the the cases where $n$ is an even integer ---
for example, in the study of knot invariants (\cite{MN}), or in logarithmic
conformal field theories (\cite{FGST1}, \cite{FGST2}).  In this paper, 
we work out a fairly detailed study on the category of finite dimensional 
modules of the restricted quantum $\overline{U}_q(\mathfrak{sl}_2)$ where 
$q$ is a $2p$-th root of unity, $p\ge2$.

\smallskip

Vertex operator algebras (VOAs) are axiomatic basis for conformal field
theories and, like other algebraic structures, have their own
representation theories. In order for a conformal field theory to make
sense on higher genus Riemann surfaces, the corresponding VOA should
satisfy certain finiteness conditions such as Zhu's $C_{2}$-finiteness
condition (\cite{Zhu}).

It is a nontrivial task to give examples of VOA which satisfy
$C_{2}$-finiteness condition --- among them are the triplet
$W$-algebras $W(p)$ ($p = 2, 3, \cdots$ ) (See \cite{FGST1}, \cite{FGST2} or 
\cite{TN} for the definition of $W(p)$). It is known that the category of
$W(p)$-modules is not semisimple and the conformal field theory
associated to $W(p)$ is so-called a logarithmic com-formal field theory;
the correlation functions may have logarithmic singularities, which are
not observed in semisimple conformal field theories.  Let us denote by
$W(p)$-${\boldsymbol{\mathrm{mod}}}$ the category of $W(p)$-modules. 
It is a braided tensor category via the fusion tensor products. 
Feigin et al. (\cite{FGST1}, \cite{FGST2}) make a new bridge between 
logarithmic conformal field theories and representation theory of 
the restricted quantum enveloping algebras. 
More precisely, they gave a following conjecture:

\begin{conj}[\cite{FGST2}]\label{conj:f}
Let $p\geq 2$ and $\overline{U}_q(\mathfrak{sl}_2)$ be the 
restricted quantum enveloping algebra 
associated to $\mathfrak{sl}_2$ at $2p$-th roots of unity.
As a braided quasitensor category,  
$W(p)$-${\boldsymbol{\mathrm{mod}}}$ is equivalent to
$\overline{U}_q(\mathfrak{sl}_2)$-${\boldsymbol{\mathrm{mod}}}$.
Here we denote by 
$\overline{U}_q(\mathfrak{sl}_2)$-${\boldsymbol{\mathrm{mod}}}$
the category of finite dimensional $\overline{U}_q(\mathfrak{sl}_2)$-modules.
\end{conj}
They also proved the conjecture for $p=2$. After the above 
conjecture, Tsuchiya and Nagatomo proved the following result.
\begin{theorem}[\cite{TN}]\label{thm:TN}
As abelian categories, these are equivalent for any $p\geq 2$. 
\end{theorem}

These works motivate our investigation of the ``quantum group-side'' of
the FGST's correspondence, in particular, as tensor categories. Our
paper is devoted to a detailed study of the tensor structure for
$\overline{U}_q(\mathfrak{sl}_2)$-${\boldsymbol{\mathrm{mod}}}$ at 
$2p$-th roots of unity with $p\geq 2$. 

\smallskip

This paper organized as follows. In Section \ref{sec:u_q}, 
the definition of $\overline{U}_q(\mathfrak{sl}_2)$ is recalled and 
the known facts about 
$\overline{U}_q(\mathfrak{sl}_2)$-${\boldsymbol{\mathrm{mod}}}$ are reviewed 
following
\cite{Sut}, \cite{X2}, \cite{CPrem}, \cite{FGST2} and \cite{Ari1}.  
Since $\overline{U}_q(\mathfrak{sl}_2)$ is a finite dimensional algebra, 
the technique of Auslander-Reiten theory allows us to completely classify
finite dimensional indecomposable
$\overline{U}_q(\mathfrak{sl}_2)$-modules. There exist $2p$ simple
modules (two of them are projective), $2p-2$ nonsimple indecomposable
projective modules, and several infinite sequences of other
indecomposable modules of semisimple length 2. Moreover
$\overline{U}_q(\mathfrak{sl}_2)$ has a tame representation type and the
Auslander-Reiten quiver of
$\overline{U}_q(\mathfrak{sl}_2)$-\textbf{mod} is determined.

In Section \ref{sec:calc} we give formulas for indecomposable decomposition of 
tensor products of arbitrary finite dimensional indecomposable 
$\overline{U}_q(\mathfrak{sl}_2)$-modules. Since 
$\overline{U}_q(\mathfrak{sl}_2)$ is a Hopf algebra, 
$\overline{U}_q(\mathfrak{sl}_2)$-${\boldsymbol{\mathrm{mod}}}$ has a 
natural tensor structure. 
Tensor product decomposition rules of simple and/or
projective modules are studied in \cite{Sut}. For computing tensor products
including other types of modules, the following general properties of
finite dimensional Hopf algebras (See Appendix A) are helpful:
\begin{enumerate}
 \item[(i)] If $\mathcal{P}$ is a projective 
$\overline{U}_q(\mathfrak{sl}_2)$-module, 
$\mathcal{Z}\otimes_k\mathcal{P}$ and $\mathcal{P}\otimes_k\mathcal{Z}$ are 
also projective for any $\overline{U}_q(\mathfrak{sl}_2)$-module 
$\mathcal{Z}$. 
\item[(ii)] All projective modules are injective. Conversely, all injective 
modules are projective.
\item[(iii)] The category of finite-dimensional 
$\overline{U}_q(\mathfrak{sl}_2)$-modules has a structure of 
a rigid tensor category. 
From the rigidity we have $\Ext_{\overline{U}_q(\mathfrak{sl}_2)}^n
(\mathcal{Z}_1\otimes_k \mathcal{Z}_2,\mathcal{Z}_3)\cong
\Ext_{\overline{U}_q(\mathfrak{sl}_2)}^n
\bigl(\mathcal{Z}_1,\mathcal{Z}_3\otimes_k D(\mathcal{Z}_2)\bigr)$ 
for arbitrary $\overline{U}_q(\mathfrak{sl}_2)$-modules 
$\mathcal{Z}_1$, $\mathcal{Z}_2$, and $\mathcal{Z}_3$, 
where $D(\mathcal{Z})$ is the standard dual of $\mathcal{Z}$. 
\end{enumerate}

\smallskip
By using the above facts, we can determine indecomposable decomposition of 
all tensor products of indecomposable $\overline{U}_q(\mathfrak{sl}_2)$-modules 
in explicit formulas. As a by-product, it is shown that  
$\overline{U}_q(\mathfrak{sl}_2)$-${\boldsymbol{\mathrm{mod}}}$ is {\it not} 
a braided tensor category if $p\ge3$. It is also proved that 
$\overline{U}_q(\mathfrak{sl}_2)$ has {\it no} universal $R$-matrices for 
$p\ge3$. Our result suggests that Conjecture \ref{conj:f} needs to be
modified; although $W(p)$-${\boldsymbol{\mathrm{mod}}}$ and
$\overline{U}_q(\mathfrak{sl}_2)$-${\boldsymbol{\mathrm{mod}}}$ are 
equivalent as abelian categories by Theorem \ref{thm:TN}, 
but their natural tensor structures do not agree with each other.

The resolution of this ``contradiction'' is a future problem. 
In the last section, we introduce a finite dimensional 
Hopf algebra $\overline{D}$ which contains 
$\overline{U}_q(\mathfrak{sl}_2)$ as a Hopf subalgebra. It is known that
$\overline{D}$ is quasi-triangular ;the explicit form of 
a universal $R$-matrix of $\overline{D}$ is given in \cite{FGST1}. 
We discuss a relationship between 
$\overline{U}_q(\mathfrak{sl}_2)$-${\boldsymbol{\mathrm{mod}}}$
and the category of finite 
dimensional representations of $\overline{D}$, and explain why 
$\overline{U}_q(\mathfrak{sl}_2)$ has no universal $R$-matrices for $p\ge3$.

\medskip
\begin{Ac}
Research of YS is supported by Grant-in-Aid for Scientific Research (C)
No. 2054009.
The authors are grateful to Professor Akishi Kato and Professor 
Akihiro Tsuchiya for valuable discussions. The authors also would like to
thank Professor Jie Xiao for valuable comments on the earlier draft. 
\end{Ac}

\section{Indecomposable modules over 
$\overline{U}_q(\mathfrak{sl}_2)$}%
\label{sec:u_q}
Throughout the paper, we work on a 
fixed algebraic closed field $k$ with characteristic zero. 
All modules considered are left modules and 
finite dimensional over $k$.

Let $p\ge 2$ be an integer and $q$ be a primitive $2p$-th root of unity.
For any integer $n$, we set
\[
 [n]=\frac{q^n-q^{-n}}{q-q^{-1}}. 
\]
Note that $[n]=[p-n]$ for any $n$. 

In this section we summarize facts about the restricted 
quantum $\mathfrak{sl}_2$, 
which one can find in \cite{Sut}, \cite{X2}, \cite{CPrem}, \cite{FGST2} and 
\cite{Ari1}. 
\subsection{The restricted quantum group %
$\boldsymbol{\overline{U}_q(\mathfrak{sl}_2)}$}
The restricted quantum group
$\overline{U}=\overline{U}_q(\mathfrak{sl}_2)$ is defined as 
an unital associative $k$-algebra with generators $E$, $F$, $K$, $K^{-1}$ and 
relations
\[
 KK^{-1}=K^{-1}K=1, \quad KEK^{-1}=q^2E, \quad KFK^{-1}=q^{-2}F, 
\]
\[
EF-FE=\frac{K-K^{-1}}{q-q^{-1}}, \quad
K^{2p}=1, \quad E^p=0, \quad F^p=0.
\]
This is a finite dimensional algebra and has a Hopf algebra structure,
where the coproduct $\Delta$, the counit $\varepsilon$, and the antipode $S$ 
are defined by
\begin{align*}
\Delta&\colon E\longmapsto E\otimes K+1\otimes E, \quad
F\longmapsto F\otimes 1+K^{-1}\otimes F, \\
&\phantom{\colon}K\longmapsto K\otimes K, \quad 
K^{-1}\longmapsto K^{-1}\otimes K^{-1}, \\
\varepsilon&\colon E\longmapsto 0,\quad F\longmapsto 0, \quad
K\longmapsto 1, \quad K^{-1}\longmapsto 1, \\
S&\colon E\longmapsto -EK^{-1}, \quad F\longmapsto -KF, \quad
K\longmapsto K^{-1}, \quad K^{-1}\longmapsto K. 
\end{align*}
The category $\overline{U}$-${\boldsymbol{\mathrm{mod}}}$ of 
finite dimensional left $\overline{U}$-modules has a structure of  
a monoidal category associated with this Hopf algebra structure on 
$\overline{U}$.
\subsection{Basic algebra}
Let $A$ be an unital associative $k$-algebra of finite dimension. 
The {\itshape basic algebra} of $A$ is defined as follows: 
Let $A=\bigoplus_{i=1}^n \mathcal{P}_i^{m_i}$ be a decomposition of $A$ into
indecomposable left ideals, 
where $\mathcal{P}_i\not\cong \mathcal{P}_j$ if $i\neq j$. 
For each $i$ take an idempotent $e_i\in A$ 
such that $Ae_i\cong \mathcal{P}_i$, 
and set $e=\sum_{i=1}^{n}e_i$. 
Then the subspace $B_A=eAe$ of $A$ has a natural $k$-algebra structure 
and is called the basic algebra of $A$. 

It is known (see \cite{ASS}, for example) that the categories of 
finite dimensional
modules over $A$ and $B_A$ are
equivalent each other by 
$B_A\text{-}{\boldsymbol{\mathrm{mod}}}\longrightarrow 
A\text{-}{\boldsymbol{\mathrm{mod}}}$; 
$\mathcal{Z}\longmapsto Ae\otimes_{B_A} \mathcal{Z}$. 

\medskip

The basic algebra $B_{\overline{U}}$ of $\overline{U}$ can be
decomposed as a direct product $B_{\overline{U}}\cong\prod_{s=0}^p B_s$ and 
one can describe each $B_s$ as follows:
\begin{itemize}
 \item $B_0\cong B_p\cong k$. 
 \item For each $s=1,\ldots,p-1$, $B_s$ is isomorphic to
the 8-dimensional algebra $B$ defined by the following quiver
\medskip
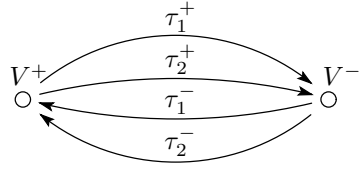
\begin{center}
 \input{quiver.tpc}
\end{center}
\medskip
with relations $\tau_i^{\pm}\tau_i^{\mp}=0$ for $i=1,2$, and
$\tau_1^{\pm}\tau_2^{\mp}=\tau_2^{\pm}\tau_1^{\mp}$.
\end{itemize}

The algebra $B$ is studied in \cite{Sut} and \cite{X2} and is known to 
have a tame representation type. 
We shall review on the classification theorem of 
isomorphism classes of indecomposable $B$-modules. 
Note that one can identify a $B$-module with 
data $\mathcal{Z}=(V_{\mathcal{Z}}^+,V_{\mathcal{Z}}^-, 
\tau_{1,\mathcal{Z}}^+,\tau_{2,\mathcal{Z}}^+,
\tau_{1,\mathcal{Z}}^-,\tau_{2,\mathcal{Z}}^-)$, 
where $V_{\mathcal{Z}}^{\pm}$ is a vector space over $k$ and 
$\tau_{i,\mathcal{Z}}^{\pm}\colon 
V_{\mathcal{Z}}^{\pm}\longrightarrow V_{\mathcal{Z}}^{\mp}$ 
($i=1,2$) are 
$k$-linear maps satisfying 
$\tau_{i,\mathcal{Z}}^{\pm}\tau_{i,\mathcal{Z}}^{\mp}=0$, 
$\tau_{1,\mathcal{Z}}^{\pm}\tau_{2,\mathcal{Z}}^{\mp}
=\tau_{2,\mathcal{Z}}^{\pm}\tau_{1,\mathcal{Z}}^{\mp}$. 

\begin{prop}\label{prop:indec B-modules}
Any indecomposable $B$-module is isomorphic to 
exactly one of modules in the following list:
\begin{itemize}
 \item Simple modules 
\[
\mathcal{X}^+=(k,0,0,0,0,0), \quad
\mathcal{X}^-=(0,k,0,0,0,0).
\]
  \item Projective-injective modules 
\[
\mathcal{P}^+=(k^2,k^2,e_{1,1},e_{2,1},e_{2,2},e_{2,1}), \quad
\mathcal{P}^-=(k^2,k^2,e_{2,2},e_{2,1},e_{1,1},e_{2,1}), 
\]
where for positive integers $m,n$ and 
$i=1,\ldots,m$, $j=1,\ldots,n$ we denote the composition of 
$j$-th projection and $i$-th embedding 
$k^n\longrightarrow k\longrightarrow k^m$ by $e_{i,j}$. 
 \item $\mathcal{M}^+(n)=\bigl(k^{n-1},k^n,\sum_{i=1}^{n-1}e_{i,i},
\sum_{i=1}^{n-1}e_{i+1,i},0,0\bigr)$, 
$\mathcal{M}^-(n)=\bigl(k^{n},k^{n-1},0,0,\sum_{i=1}^{n-1}e_{i,i},
\sum_{i=1}^{n-1}e_{i+1,i}\bigr)$ 
for each integer $n\ge2$. 
 \item $\mathcal{W}^+(n)=\bigl(k^n,k^{n-1},\sum_{i=1}^{n-1}e_{i,i},
\sum_{i=1}^{n-1}e_{i,i+1},0,0\bigr)$, 
$\mathcal{W}^-(n)=\bigl(k^{n-1},k^n,0,0,\sum_{i=1}^{n-1}e_{i,i},
\sum_{i=1}^{n-1}e_{i,i+1}\bigr)$ 
for each integer $n\ge2$. 
 \item $\mathcal{E}^+(n;\lambda)=\bigl(k^n,k^n,\varphi_1(n;\lambda),
\varphi_2(n;\lambda),0,0\bigr)$, 
$\mathcal{E}^-(n;\lambda)=\bigl(k^n,k^n,0,0,\varphi_1(n;\lambda),
\varphi_2(n;\lambda)\bigr)$  
for each integer $n\ge1$ and $\lambda\in\mathbb{P}^1(k)$, 
where 
\[
\bigl(\varphi_1(n;\lambda),\varphi_2(n;\lambda)\bigr)
=\begin{cases}
\bigl(\beta\cdot\id+\sum_{i=1}^{n-1}e_{i,i+1},\id\bigr)
&(\lambda=[\beta:1]),\\
\bigl(\id,\sum_{i=1}^{n-1}e_{i,i+1}\bigr)&(\lambda=[1:0]).
\end{cases}
\]
\end{itemize}
\end{prop}

\subsection{Indecomposable modules}
\begin{defn}
For $s=1,\ldots,p-1$, Let $\Phi_s$ be the composition of functors
$B$-${\boldsymbol{\mathrm{mod}}}\longrightarrow
B_{\overline{U}}$-${\boldsymbol{\mathrm{mod}}}\longrightarrow
\overline{U}$-${\boldsymbol{\mathrm{mod}}}$, 
where the first one is induced from 
$B_{\overline{U}}\cong\prod_{s=0}^pB_s\longrightarrow B_s\cong B$  
and the second one is expressed in the previous subsection. 

We denote by $\mathcal{X}_s^+$, $\mathcal{X}_{p-s}^-$, 
$\mathcal{P}_s^+$, $\mathcal{P}_{p-s}^-$, $\mathcal{M}_s^+(n)$, 
$\mathcal{M}_{p-s}^-(n)$, 
$\mathcal{W}_s^+(n)$, $\mathcal{W}_{p-s}^-(n)$, 
$\mathcal{E}_s^+(n;\lambda)$, $\mathcal{E}_{p-s}^-(n;\lambda)$
the images of $\mathcal{X}^+$, $\mathcal{X}^-$, 
$\mathcal{P}^+$, $\mathcal{P}^-$, $\mathcal{M}^+(n)$, $\mathcal{M}^-(n)$, 
$\mathcal{W}^+(n)$, $\mathcal{W}^-(n)$, 
$\mathcal{E}^+(n;\lambda)$, $\mathcal{E}^-(n;\lambda)$ by $\Phi_s$. 
\end{defn}

Denote by $\mathcal{C}(s)$ the full subcategory of 
$\overline{U}$-${\boldsymbol{\mathrm{mod}}}$ corresponding to 
$B_s$-modules (considered as $B_{\overline{U}}$-modules) 
for $s=0,\ldots,p$.
Each indecomposable $\overline{U}$-module belongs to exactly one of 
$\mathcal{C}(s)$ ($s=0,\ldots, p$).

Since $B_0\cong B_p\cong k$, each of $\mathcal{C}(0)$ and $\mathcal{C}(p)$ has 
precisely one indecomposable 
module (denoted by $\mathcal{X}_p^+$, $\mathcal{X}_p^-$, respectively).

For $s=1,\ldots,p-1$, indecomposable modules in $\mathcal{C}(s)$ are 
classified as follows.

\begin{prop}
Each subcategory $\mathcal{C}(s)$ $(s=1,\ldots,p-1)$ has two 
simple modules $\mathcal{X}_s^+$ and $\mathcal{X}_{p-s}^-$, 
two indecomposable projective-injective modules 
$\mathcal{P}_s^+$ and $\mathcal{P}_{p-s}^-$, 
and three series of indecomposable modules:
\begin{itemize}
 \item $\mathcal{M}_s^+(n)$ and $\mathcal{M}_{p-s}^-(n)$ for each
integer $n\ge2$, 
 \item $\mathcal{W}_s^+(n)$ and $\mathcal{W}_{p-s}^-(n)$ for each
integer $n\ge2$, 
 \item $\mathcal{E}_s^+(n;\lambda)$ and $\mathcal{E}_{p-s}^-(n;\lambda)$ for each
integer $n\ge1$ and $\lambda\in\mathbb{P}^1(k)$, 
\end{itemize}
Moreover any indecomposable module in $\mathcal{C}(s)$ is isomorphic to 
one of the modules listed above.
\end{prop}

Since a complete set of primitive orthogonal idempotents of $\overline{U}$ is
known (see \cite{Ari1}, for example), we can describe all the above 
indecomposable modules explicitly by bases and action of $\overline{U}$ on 
those. However, we give them only for $\mathcal{X}_s^{\pm}$ $(s=1,\ldots,p)$
and $\mathcal{E}_s^{\pm}(1;\lambda)$ 
$(s=1,\ldots,p-1,\,\lambda=[\lambda_1:\lambda_2]\in \mathbb{P}^1(k))$ in the
next proposition, 
because it is enough for computing tensor products of indecomposable modules.
\begin{prop}\label{prop:basis}
{\rm (i)}\
$\mathcal{X}_s^{\pm}$ $(s=1,\ldots,p)$ 
is isomorphic to the $s$-dimensional module defined by 
basis $\{a_n\}_{n=0,\ldots,s-1}$ and $\overline{U}$-action given by
\[
Ka_n=\pm q^{s-1-2n}a_n, \quad
Ea_n=\begin{cases}\pm[n][s-n]a_{n-1}&(n\neq0)\\0&(n=0)\end{cases},\quad
Fa_n=\begin{cases}a_{n+1}&(n\neq s-1)\\0&(n=s-1)\end{cases}.
\]
{\rm (ii)}\ $\mathcal{E}_s^{\pm}(1;\lambda)$ 
$(s=1,\ldots,p-1,\,\lambda=[\lambda_1:\lambda_2])$ is isomorphic to the 
$p$-dimensional module defined by 
basis $\{b_n\}_{n=0,\ldots,s-1}\amalg\{x_m\}_{m=0,\ldots,p-s-1}$ and 
$\overline{U}$-action given by
\begin{align*}
Kb_n&=\pm q^{s-1-2n}b_n, \quad Kx_m=\mp q^{p-s-1-2m}x_m, \\
Eb_n&=\begin{cases}\pm[n][s-n]b_{n-1}&(n\neq0)\\
\lambda_2x_{p-s-1}&(n=0)\end{cases},\quad
Ex_m=\begin{cases}\mp[m][p-s-m]x_{m-1}&(m\neq0)\\0&(m=0)\end{cases},\\
Fb_n&=\begin{cases}b_{n+1}&(n\neq s-1)\\
\lambda_1x_0&(n=s-1)\end{cases}, \quad
Fx_m=\begin{cases}x_{m+1}&(m\neq p-s-1)\\0&(m=p-s-1)\end{cases}.
\end{align*}
\end{prop} 

We shall introduce some basic notations in representation theory
of finite dimensional algebras.
\begin{defn}
Let $A$ be a unital associative $k$-algebra of finite dimension and 
$\mathcal{Z}$ a finite dimensional left $A$-module.
\vskip 1mm
\noindent{\rm (i)}\ 
The radical $\mbox{\rm rad}\mathcal{Z}$ of $\mathcal{Z}$ is the
intersection of all the maximal proper submodules of $\mathcal{Z}$.

\noindent{\rm (ii)}\ The module $\mathcal{Z}/\mbox{\rm rad}\mathcal{Z}$ is 
the largest semisimple factor module of $\mathcal{Z}$ which is called the 
top of $\mathcal{Z}$. We denote it $\mbox{\rm top}\mathcal{Z}$.

\noindent{\rm (iii)}\ The sum of all simple submodules of $\mathcal{Z}$ is
called the socle of $\mathcal{Z}$ which is denoted by $\soc \mathcal{Z}$. 

\noindent{\rm (iv)}\ We define a semisimple filtration of $\mathcal{Z}$ as
a sequence of submodules
$$\mathcal{Z}=\mathcal{Z}_0\supset \mathcal{Z}_1\supset \cdots \supset
\mathcal{Z}_l=0$$
such that each quotient $\mathcal{Z}_i/\mathcal{Z}_{i+1}$ is semisimple.
The number $l$ is called the length of the filtration. In the set of 
semisimple filtrations of $\mathcal{Z}$, there exists a filtration with the
minimum length $l$. We call $l$ the semisimple length of $\mathcal{Z}$.
We remark that an indecomposable module with semisimple length $1$ is nothing
but a simple module.
\end{defn}

Let us return to our case.
\begin{prop}
{\rm (i)}\
There are no $\overline{U}$-modules with semisimple length greater than $3$.

\noindent{\rm (ii)}\ 
The only indecomposable modules with semisimple length $3$ are
the projective modules $\mathcal{P}_s^{\pm}$ with $s=1,\ldots,p-1$. 
More precisely, for $s=1,\ldots,p-1$, the projective module 
$\mathcal{P}_s^{\pm}$ has the following semisimple filtration with length $3$:
$$\mathcal{P}_s^{\pm}=(\mathcal{P}_s^{\pm})_0\supset
(\mathcal{P}_s^{\pm})_1\supset (\mathcal{P}_s^{\pm})_2\supset 
(\mathcal{P}_s^{\pm})_3=0$$
such that 
$$(\mathcal{P}_s^{\pm})_0/(\mathcal{P}_s^{\pm})_1=
\top \mathcal{P}_s^{\pm}\cong\mathcal{X}_{s}^{\pm},\quad
(\mathcal{P}_s^{\pm})_1/(\mathcal{P}_s^{\pm})_2\cong 
(\mathcal{X}_{p-s}^{\mp})^{2},\quad
(\mathcal{P}_s^{\pm})_2=\soc \mathcal{P}_s^{\pm}\cong\mathcal{X}_{s}^{\pm}.$$

\noindent{\rm (iii)}\ The other non-simple indecomposable modules have
semisimple length $2$. More precisely, for $s=1,\ldots,p-1$, we have 
$$\top \mathcal{M}_s^{\pm}(n)\cong(\mathcal{X}_{s}^{\pm})^{n-1},\quad  
\top \mathcal{W}_s^{\pm}(n)\cong(\mathcal{X}_{s}^{\pm})^n,\quad 
\top \mathcal{E}_s^{\pm}(n;\lambda)\cong(\mathcal{X}_{s}^{\pm})^n,$$  
$$\soc \mathcal{M}_s^{\pm}(n)\cong(\mathcal{X}_{p-s}^{\mp})^n,\quad 
\soc \mathcal{W}_s^{\pm}(n)\cong(\mathcal{X}_{p-s}^{\mp})^{n-1},\quad  
\soc \mathcal{E}_s^{\pm}(n;\lambda)\cong(\mathcal{X}_{p-s}^{\mp})^n.$$  
\end{prop}

\begin{cor} 
We have
$\dim_k \mathcal{X}_s^{\pm}=s$, 
$\dim_k \mathcal{P}_s^{\pm}=2p$, 
$\dim_k \mathcal{M}_s^{\pm}(n)=pn-s$, 
$\dim_k \mathcal{W}_s^{\pm}(n)=pn-p+s$, 
$\dim_k \mathcal{E}_s^{\pm}(n;\lambda)=pn$.  
\end{cor}

\subsection{Extensions}
We describe the projective covers and the injective envelopes of 
indecomposable $\overline{U}$-modules 
which we use in the sequel. 

\begin{prop}\label{prop:covers}
There exist following exact sequences
\[
 0\longrightarrow \mathcal{M}_{p-s}^{\mp}(n)
\longrightarrow (\mathcal{P}_{s}^{\pm})^n
\longrightarrow \mathcal{M}_{s}^{\pm}(n+1)
\longrightarrow 0,
\]
\[
 0\longrightarrow \mathcal{W}_{p-s}^{\mp}(n+1)
\longrightarrow (\mathcal{P}_{s}^{\pm})^n
\longrightarrow \mathcal{W}_{s}^{\pm}(n)
\longrightarrow 0,
\]
\[
 0\longrightarrow \mathcal{E}_{p-s}^{\mp}(n;-\lambda)
\longrightarrow (\mathcal{P}_{s}^{\pm})^n
\longrightarrow \mathcal{E}_{s}^{\pm}(n;\lambda)
\longrightarrow 0
\]
for each $s=1,\ldots,p-1$, $n\ge1$ and $\lambda\in\mathbb{P}^1(k)$, 
where we set $\mathcal{M}_{p-s}^{\mp}(1)=
\mathcal{W}_{s}^{\pm}(1)=\mathcal{X}_s^{\pm}$. 
Moreover, each sequence gives 
the projective cover of the right term and 
the injective envelope of the left term.
\end{prop}

The first extensions between indecomposable $\overline{U}$-modules 
can be calculated by passing to $B$-${\boldsymbol{\mathrm{mod}}}$ and 
using the Auslander-Reiten formulas (\cite{ASS}).

\begin{prop}\label{prop:ext}
{\rm (i)}\
$\Ext_{\overline{U}}^1
\bigl(\mathcal{E}_s^{\pm}(n;\lambda),\mathcal{X}_s^{\pm}\bigr)=0$, 
$\dim_k\Ext_{\overline{U}}^1
\bigl(\mathcal{E}_s^{\pm}(n;\lambda),\mathcal{X}_{p-s}^{\mp}\bigr)=n$.

\noindent{\rm (ii)}\
$\dim_k\Ext_{\overline{U}}^1
\bigl(\mathcal{X}_s^{\pm},\mathcal{E}_s^{\pm}(n;\lambda)\bigr)=n$, 
$\Ext_{\overline{U}}^1
\bigl(\mathcal{X}_{p-s}^{\mp},\mathcal{E}_s^{\pm}(n;\lambda)\bigr)=0$.

\noindent{\rm (iii)}\ 
$\dim_k\Ext_{\overline{U}}^1
\bigl(\mathcal{E}_s^{\pm}(m;\lambda),\mathcal{E}_{s}^{\pm}(n;\mu)\bigr)
=\delta_{\lambda\mu}\min\{m,n\}$, 
$\dim_k\Ext_{\overline{U}}^1
\bigl(\mathcal{E}_s^{\pm}(m;\lambda),\mathcal{E}_{p-s}^{\mp}(n;-\mu)\bigr)
=\delta_{\lambda\mu}\min\{m,n\}$.  
\end{prop}
For later use, the following exact sequences are also useful.
\begin{prop}\label{prop:exact seq for E}
Let $s=1,\ldots,p-1$, $n\ge2$ and $\lambda\in\mathbb{P}^1(k)$. 
%
%
Then there exist exact sequences
\[
 0\longrightarrow \mathcal{E}_{s}^{\pm}(n-1;\lambda)
\longrightarrow \mathcal{E}_{s}^{\pm}(n;\lambda)
\longrightarrow \mathcal{E}_{s}^{\pm}(1;\lambda)
\longrightarrow 0. 
\]
\end{prop}
\section{Calculation of tensor products}\label{sec:calc}
\subsection{Tensor products of simple modules}
Tensor products of simple $\overline{U}$-modules 
$\mathcal{X}_s^{\pm}\otimes\mathcal{X}_{s'}^{\pm}$ 
($-\otimes-$ means $-\otimes_k-$, here and further)
have been studied in \cite{Sut}. 
Here we present these results with some different notation.

\begin{defn}
For $s, s'=1,\ldots, p$ with $s\le s'$, define $I_{s,s'}$ and
$J_{s,s'}$ by
\begin{align*}
I_{s,s'}&=\{t=s'-s+2i-1\mid i=1,\ldots, s,\ t\le 2p-s-s'\}, \\
J_{s,s'}&=\{t=2p-2i-s'+s+1\mid i=1,\ldots, s,\ t\le p\}, 
\end{align*}
and set $I_{s,s'}=I_{s',s}$, $J_{s,s'}=J_{s',s}$ for 
$s, s'=1,\ldots, p$ with $s>s'$. 
\end{defn}

\begin{ex} Let $p=5$. Then $I_{s,s'}$ and $J_{s,s'}$ are
as the following table.

\medskip
\begin{center}
 \begin{tabular}{c|ccccc}
 $I$ & $1$ & $2$ & $3$ & $4$ & $5$ \\\hline
 $1$ & $\{1\}$ & $\{2\}$ & $\{3\}$ & $\{4\}$ & $\emptyset$\\
 $2$ & $\{2\}$ & $\{1,3\}$ & $\{2,4\}$ & $\{3\}$ & $\emptyset$\\
 $3$ & $\{3\}$ & $\{2,4\}$ & $\{1,3\}$ & $\{2\}$ & $\emptyset$\\
 $4$ & $\{4\}$ & $\{3\}$ & $\{2\}$ & $\{1\}$ & $\emptyset$\\
 $5$ & $\emptyset$ & $\emptyset$ & $\emptyset$ & $\emptyset$ & $\emptyset$
 \end{tabular}\qquad
 \begin{tabular}{c|ccccc}
 $J$ & $1$ & $2$ & $3$ & $4$ & $5$ \\\hline
 $1$ & $\emptyset$ & $\emptyset$ & $\emptyset$ & $\emptyset$ & $\{5\}$\\
 $2$ & $\emptyset$ & $\emptyset$ & $\emptyset$ & $\{5\}$ & $\{4\}$\\
 $3$ & $\emptyset$ & $\emptyset$ & $\{5\}$ & $\{4\}$ & $\{3,5\}$\\
 $4$ & $\emptyset$ & $\{5\}$ & $\{4\}$ & $\{3,5\}$ & $\{2,4\}$\\
 $5$ & $\{5\}$ & $\{4\}$ & $\{3,5\}$ & $\{2,4\}$ & $\{1,3,5\}$
 \end{tabular}
\end{center}
\end{ex}
\medskip

We collect some properties of $I_{s,s'}$ and $J_{s,s'}$ for 
later use, a proof of which is straightforward. 

\begin{prop}\label{prop:I and J}
Let $s, s', t, t'=1,\ldots,p$. 

\noindent{\rm (i)}\ 
$I_{s,s'}\subset \{1,\ldots,p-1\}$, $J_{s,s'}\subset \{1,\ldots,p\}$.

\noindent{\rm (ii)}\ 
$I_{s,s'}\cap J_{s,s'}=\emptyset$.

\noindent{\rm (iii)}\
If $s=1,\ldots, p-1$, $I_{p-s,s'}=\{p-t\mid t\in I_{s,s'}\}$. 
If $s=p$, $I_{p,s'}=\emptyset$. 

\noindent{\rm (iv)}\
$t\in I_{s,s'}$ implies $s'\in I_{s,t}$. 

\noindent{\rm (v)}\
$J_{s,s'}=J_{t,t'}$ if $s+s'=t+t'$. 
If $s+s'\le p$, $J_{s,s'}=\emptyset$.
\end{prop}

\begin{rem}
Since $J_{s,s'}$ depends only on $s+s'$ by (v), 
we denote it by $J_{s+s'}$ in the following. 
\end{rem}

\begin{thm}[\cite{Sut}]\label{thm:simple and simple}
For $s,s'=1,\ldots,p$ we have
\begin{align*}
\mathcal{X}_s^+\otimes\mathcal{X}_{s'}^+&\cong
\bigoplus_{t\in I_{s,s'}}\mathcal{X}_t^+\oplus
\bigoplus_{t\in J_{s+s'}}\mathcal{P}_t^+,  \\
\mathcal{X}_s^{\pm}\otimes \mathcal{X}_1^-
&\cong \mathcal{X}_1^-\otimes \mathcal{X}_s^{\pm}
\cong \mathcal{X}_s^{\mp},  \\
\mathcal{P}_s^{\pm}\otimes \mathcal{X}_1^-
&\cong \mathcal{X}_1^-\otimes \mathcal{P}_s^{\pm}
\cong \mathcal{P}_s^{\mp},
\end{align*}
where we set $\mathcal{P}_p^{\pm}=\mathcal{X}_p^{\pm}$. 
\end{thm}

\begin{rem}
The second and third formulas of the theorem enable us to
compute the tensor products 
$\mathcal{X}_s^-\otimes \mathcal{X}_{s'}^+$, 
$\mathcal{X}_s^+\otimes \mathcal{X}_{s'}^-$ and 
$\mathcal{X}_s^-\otimes \mathcal{X}_{s'}^-$. 
For example, 
$\mathcal{X}_s^-\otimes \mathcal{X}_{s'}^+
\cong \mathcal{X}_1^-\otimes \mathcal{X}_s^+\otimes \mathcal{X}_{s'}^+
\cong \mathcal{X}_1^-\otimes\bigl(\bigoplus_{t\in I_{s,s'}}\mathcal{X}_t^+\oplus
\bigoplus_{t\in J_{s+s'}}\mathcal{P}_t^+\bigr)
\cong \bigoplus_{t\in I_{s,s'}}\mathcal{X}_t^-\oplus
\bigoplus_{t\in J_{s+s'}}\mathcal{P}_t^-$. 
In the following this kind of procedure will be omitted. 
\end{rem}
\subsection{Tensor products with projective modules}\label{subsec:proj}
The tensor products of projective modules with simple modules
are also computed in \cite{Sut}:

\begin{thm}[\cite{Sut}]\label{thm:projective and simple}
For $s=1,\ldots,p-1$ and $s'=1,\ldots,p$ we have
\[
 \mathcal{P}_s^+\otimes \mathcal{X}_{s'}^+
\cong \mathcal{X}_{s'}^+\otimes \mathcal{P}_{s}^+
\cong \bigoplus_{t\in I_{s,s'}}\mathcal{P}_t^+\oplus
\bigoplus_{t\in J_{s+s'}}(\mathcal{P}_t^+)^2\oplus
\bigoplus_{t\in J_{p-s+s'}}(\mathcal{P}_t^-)^2.
\]
\end{thm}

Let us calculate the tensor products of projective modules with
arbitrary modules. 

\begin{cor}\label{cor:projective and arbitrary}
Suppose $s=1,\ldots, p-1$. Let $\mathcal{Z}$ be an arbitrary 
$\overline{U}$-module and $\oplus_{i\in\Lambda}\mathcal{S}_i$ the 
the direct sum of its composition factors of $\mathcal{Z}$. Then we have

\noindent {\rm (i)}\
$\mathcal{P}_s^{\pm}\otimes \mathcal{Z}\cong 
\oplus_{i\in\Lambda}\mathcal{P}_s^{\pm}\otimes\mathcal{S}_i$ and 
$\mathcal{Z}\otimes \mathcal{P}_s^{\pm}\cong
\oplus_{i\in\Lambda}\mathcal{S}_i\otimes\mathcal{P}_s^{\pm}$,

\noindent{\rm (ii)}\
$\mathcal{P}_s^{\pm}\otimes \mathcal{Z}\cong\mathcal{Z}\otimes
\mathcal{P}_s^{\pm}$.
\end{cor}

\begin{proof}
The statement (i) is a direct consequence of Corollary \ref{cor:proj} in 
Appendix A. Therefore, for showing (ii), it is enough to prove that 
$\mathcal{P}_s^{\pm}\otimes \mathcal{S}\cong
\mathcal{S}\otimes\mathcal{P}_s^{\pm}$ for each
simple module $\mathcal{S}$. However, it is already proved in Theorem
\ref{thm:projective and simple}. 
\end{proof}

\begin{ex}
For $s,s'=1,\ldots,p-1$ and $n\ge2$ we have
\begin{align*}
&\mathcal{P}_s^+\otimes\mathcal{M}_{s'}^+(n)\cong
\mathcal{M}_{s'}^+(n)\otimes\mathcal{P}_s^+\\&\qquad\cong
\mathcal{P}_s^+\otimes
\bigl((\mathcal{X}_{p-s'}^-)^n\oplus(\mathcal{X}_{s'}^+)^{n-1}\bigr)\\
&\qquad\cong
\bigoplus_{t\in I_{s,s'}}
\bigl((\mathcal{P}_t^+)^{n-1}\oplus (\mathcal{P}_{p-t}^-)^{n}\bigr)\oplus
\bigoplus_{t\in J_{s+s'}}(\mathcal{P}_t^+)^{2n-2}\oplus
\bigoplus_{t\in J_{2p-s-s'}}(\mathcal{P}_t^+)^{2n}\\
&\qquad\qquad\qquad\oplus
\bigoplus_{t\in J_{p+s-s'}}(\mathcal{P}_t^-)^{2n}\oplus
\bigoplus_{t\in J_{p-s+s'}}(\mathcal{P}_t^-)^{2n-2}, 
\end{align*}
where in the last isomorphism we use Proposition \ref{prop:I and J} (iii).
\end{ex}

\subsection{Tensor products with %
$\boldsymbol{\mathcal{M}_s^{\pm}(n)}$ and %
$\boldsymbol{\mathcal{W}_s^{\pm}(n)}$}
Define a multiplicative law $\cdot$ on the set $\{+,-\}$ by
$$+\cdot +=+,\quad +\cdot-=-,\quad -\cdot +=-,\quad -\cdot -=+.$$
Namely, we regard the set $\{+,-\}$ with the multiplicative law $\cdot$ 
as $\mathbb{Z}/2\mathbb{Z}$.

\begin{thm}\label{thm:M and M}
Assume $s,s'=1,\ldots,p-1$ and $m,n\ge2$. Let $\alpha,\beta\in\{+,-\}$.
Then we have 
\begin{align*}
&\mathcal{M}_s^{\alpha}(n)\otimes\mathcal{X}_{s'}^{\beta}
\cong \mathcal{X}_{s'}^{\beta}\otimes\mathcal{M}_s^{\alpha}(n)
\cong
\bigoplus_{t\in I_{s,s'}}\mathcal{M}_{t}^{\alpha\cdot\beta}(n)
\oplus\bigoplus_{t\in J_{s+s'}}(\mathcal{P}_t^{\alpha\cdot\beta})^{n-1}
\oplus\bigoplus_{t\in J_{p-s+s'}}(\mathcal{P}_t^{-\alpha\cdot\beta})^{n}, \\
&\mathcal{W}_s^{\alpha}(n)\otimes\mathcal{X}_{s'}^{\beta}
\cong \mathcal{X}_{s'}^{\beta}\otimes\mathcal{W}_s^{\alpha}(n)
\cong 
\bigoplus_{t\in I_{s,s'}}\mathcal{W}_{t}^{\alpha\cdot\beta}(n)
\oplus\bigoplus_{t\in J_{s+s'}}(\mathcal{P}_t^{\alpha\cdot\beta})^{n}
\oplus\bigoplus_{t\in J_{p-s+s'}}(\mathcal{P}_t^{-\alpha\cdot\beta})^{n-1}, \\
&\mathcal{M}_s^{\alpha}(m)\otimes\mathcal{M}_{s'}^{\beta}(n)\\
&\qquad\cong
\bigoplus_{t\in I_{s,s'}}
\bigl(\mathcal{M}_{p-t}^{-\alpha\cdot\beta}(m+n-1)\oplus 
(\mathcal{P}_{t}^{\alpha\cdot\beta})^{(m-1)(n-1)}\bigr)
\oplus\bigoplus_{t\in J_{s+s'}}(\mathcal{P}_t^{\alpha\cdot\beta})^{(m-1)(n-1)}
\oplus\bigoplus_{t\in J_{2p-s-s'}}(\mathcal{P}_t^{\alpha\cdot\beta})^{mn}\\
&\qquad\qquad\qquad\oplus
\bigoplus_{t\in J_{p+s-s'}}(\mathcal{P}_t^{-\alpha\cdot\beta})^{(m-1)n}\oplus
\bigoplus_{t\in J_{p-s+s'}}(\mathcal{P}_t^{-\alpha\cdot\beta})^{m(n-1)}, \\
&\mathcal{W}_s^{\alpha}(m)\otimes\mathcal{W}_{s'}^{\beta}(n)\\
&\qquad\cong
\bigoplus_{t\in I_{s,s'}}
\bigl(\mathcal{W}_{t}^{\alpha\cdot\beta}(m+n-1)\oplus 
(\mathcal{P}_{t}^{\alpha\cdot\beta})^{(m-1)(n-1)}\bigr)
\oplus\bigoplus_{t\in J_{s+s'}}(\mathcal{P}_t^{\alpha\cdot\beta})^{mn}
\oplus\bigoplus_{t\in J_{2p-s-s'}}(\mathcal{P}_t^{\alpha\cdot\beta})^{(m-1)(n-1)}\\
&\qquad\qquad\qquad\oplus
\bigoplus_{t\in J_{p+s-s'}}(\mathcal{P}_t^{-\alpha\cdot\beta})^{m(n-1)}\oplus
\bigoplus_{t\in J_{p-s+s'}}(\mathcal{P}_t^{-\alpha\cdot\beta})^{(m-1)n}, \\
&\mathcal{M}_s^{\alpha}(m)\otimes\mathcal{W}_{s'}^{\beta}(n)\cong
\mathcal{W}_{s'}^{\beta}(m)\otimes\mathcal{M}_s^{\alpha}(n)\\
&\qquad\cong
\bigoplus_{t\in I_{s,s'}}\mathcal{Y}_t^{\alpha\cdot\beta}(m,n)
\oplus\bigoplus_{t\in J_{s+s'}}(\mathcal{P}_t^{\alpha\cdot\beta})^{(m-1)n}
\oplus\bigoplus_{t\in J_{2p-s-s'}}(\mathcal{P}_t^{\alpha\cdot\beta})^{m(n-1)}\\
&\qquad\qquad\qquad\oplus
\bigoplus_{t\in J_{p+s-s'}}(\mathcal{P}_t^{-\alpha\cdot\beta})^{(m-1)(n-1)}\oplus
\bigoplus_{t\in J_{p-s+s'}}(\mathcal{P}_t^{-\alpha\cdot\beta})^{mn}, 
\end{align*}
where $\mathcal{Y}_t^{\alpha}(m,n)$ is defined by
\[
 \mathcal{Y}_t^{\alpha}(m,n)=\begin{cases}
\mathcal{M}_t^{\alpha}(m-n+1)\oplus (\mathcal{P}_t^{\alpha})^{m(n-1)} 
& \text{if $m>n$},\\
\mathcal{X}_{p-t}^{-\alpha}\oplus (\mathcal{P}_t^{\alpha})^{n(n-1)} 
& \text{if $m=n$},\\
\mathcal{W}_{p-t}^{-\alpha}(n-m+1)\oplus (\mathcal{P}_t^{\alpha})^{(m-1)n} 
& \text{if $m<n$}.
\end{cases}
\]
Here we set $\mathcal{M}_{p-s}^{\mp}(1)=
\mathcal{W}_{s}^{\pm}(1)=\mathcal{X}_s^{\pm}$ as before. 
\end{thm}
\begin{proof}
We only prove the first formula. The others are proved by similar method.

Since $\mathcal{M}_{s}^{\pm}(1)=\mathcal{X}_{p-s}^{\mp}$,  
the formula is already given in Theorem \ref{thm:simple and simple} for $n=1$.
Suppose that the formula holds for $n-1$. Applying the exact functor 
$-\otimes \mathcal{X}_{s'}^{\beta}$ to 
the first exact sequence in Proposition \ref{prop:covers}, we have
an exact sequence
\[
 0\longrightarrow 
\mathcal{M}_{p-s}^{-\alpha}(n-1)\otimes\mathcal{X}_{s'}^{\beta}\longrightarrow
(\mathcal{P}_{s}^{\alpha})^{n-1}\otimes\mathcal{X}_{s'}^{\beta}\longrightarrow
\mathcal{M}_{s}^{\alpha}(n)\otimes\mathcal{X}_{s'}^{\beta}\longrightarrow 0.
\]
By the hypothesis and Proposition \ref{prop:I and J} (iii) we have
\begin{align*}
\mathcal{M}_{p-s}^{-\alpha}(n-1)\otimes\mathcal{X}_{s'}^{\beta}&\cong
\bigoplus_{t\in I_{p-s,s'}}\mathcal{M}_{t}^{-\alpha\cdot\beta}(n-1)
\oplus\bigoplus_{t\in J_{p-s+s'}}(\mathcal{P}_t^{-\alpha\cdot\beta})^{n-2}
\oplus\bigoplus_{t\in J_{s+s'}}(\mathcal{P}_t^{\alpha\cdot\beta})^{n-1}\\
&\cong
\bigoplus_{t\in I_{s,s'}}\mathcal{M}_{p-t}^{-\alpha\cdot\beta}(n-1)
\oplus\bigoplus_{t\in J_{p-s+s'}}(\mathcal{P}_t^{-\alpha\cdot\beta})^{n-2}
\oplus\bigoplus_{t\in J_{s+s'}}(\mathcal{P}_t^{\alpha\cdot\beta})^{n-1}
\end{align*}
On the other hand, we can calculate the middle term by Theorem 
\ref{thm:projective and simple}:
\begin{align*}
(\mathcal{P}_{s}^{\alpha})^{n-1}\otimes\mathcal{X}_{s'}^{\beta}&\cong
\left(\mathcal{P}_{s}^{\alpha} \otimes\mathcal{X}_{s'}^{\beta}\right)^{n-1}\\
&\cong
\bigoplus_{t\in I_{s,s'}}(\mathcal{P}_t^{\alpha\cdot\beta})^{n-1}\oplus
\bigoplus_{t\in J_{s+s'}}(\mathcal{P}_t^{\alpha\cdot\beta})^{2(n-1)}\oplus
\bigoplus_{t\in J_{p-s+s'}}(\mathcal{P}_t^{-\alpha\cdot\beta})^{2(n-1)}.
\end{align*} 
Since all projective modules are injective, the projective summands in the 
left term also appear in the middle term. Therefore we have an exact sequence
\begin{align*}
&0 \longrightarrow \bigoplus_{t\in I_{s,s'}}
\mathcal{M}_{p-t}^{-\alpha\cdot\beta}(n-1)
\longrightarrow \bigoplus_{t\in I_{s,s'}}
(\mathcal{P}_t^{\alpha\cdot\beta})^{n-1}\oplus
\bigoplus_{t\in J_{s+s'}}(\mathcal{P}_t^{\alpha\cdot\beta})^{(n-1)}\oplus
\bigoplus_{t\in J_{p-s+s'}}(\mathcal{P}_t^{-\alpha\cdot\beta})^{n}\\
&\hspace*{100mm}
\longrightarrow \mathcal{M}_{s}^{\alpha}(n)\otimes\mathcal{X}_{s'}^{\beta}
\longrightarrow 0.
\end{align*}
An injective homomorphism from 
$\mathcal{M}_{p-t}^{-\alpha\cdot\beta}(n-1)$ to an injective module must factor 
through its injective envelope $(\mathcal{P}_t^{\alpha\cdot\beta})^{n-1}$.
Consequently we have
\begin{align*}
\mathcal{M}_{s}^{\alpha}(n)\otimes\mathcal{X}_{s'}^{\beta}&\cong
\bigoplus_{t\in I_{s,s'}}\left((\mathcal{P}_t^{\alpha\cdot\beta})^{n-1}/
\mathcal{M}_{p-t}^{-\alpha\cdot\beta}(n-1)\right)
\oplus\bigoplus_{t\in J_{s+s'}}(\mathcal{P}_t^{\alpha\cdot\beta})^{n-1}
\oplus\bigoplus_{t\in J_{p-s+s'}}(\mathcal{P}_t^{-\alpha\cdot\beta})^{n}\\
&\cong
\bigoplus_{t\in I_{s,s'}}\mathcal{M}_{t}^{\alpha\cdot\beta}(n)
\oplus\bigoplus_{t\in J_{s+s'}}(\mathcal{P}_t^{\alpha\cdot\beta})^{n-1}
\oplus\bigoplus_{t\in J_{p-s+s'}}(\mathcal{P}_t^{-\alpha\cdot\beta})^{n}.
\end{align*}
For the case of $\mathcal{X}_{s'}^{\beta}\otimes\mathcal{M}_{s}^{\alpha}(n)$, 
we can determine the decomposition rule by the similar method.  
\end{proof}
\subsection{Tensor products of %
$\boldsymbol{\mathcal{E}_s^{\pm}(1;\lambda)}$ with %
simple modules}
The aim of this subsection is to compute the decomposition of
$\mathcal{E}_s^{\pm}(1;\lambda)\otimes\mathcal{X}_{s'}^{\pm}$ and 
$\mathcal{X}_{s'}^{\pm}\otimes\mathcal{E}_s^{\pm}(1;\lambda)$.
Firstly, we shall calculate tensor product of $\mathcal{E}_s^{\pm}(1;\lambda)$
and $1$-dimensional module. Let us introduce a map $\kappa:\{+,-\}\to 
\{\pm 1\}$ by
$$\kappa(+)=1\quad\mbox{and}\quad \kappa(-)=-1.$$

\begin{prop}\label{prop:E and X_1^-}
Let $\alpha,\beta\in\{+,-\}$. For $s=1,\ldots,p-1$ and 
$\lambda\in\mathbb{P}^1(k)$ we have
\begin{align*}
\mathcal{E}_s^{\alpha}(1;\lambda)\otimes\mathcal{X}_1^{\beta}&\cong
\mathcal{E}_s^{\alpha\cdot\beta}(1;\kappa(\beta)\lambda), \\
\mathcal{X}_1^{\beta}\otimes\mathcal{E}_s^{\alpha}(1;\lambda)&\cong
\mathcal{E}_s^{\alpha\cdot\beta}\bigl(1;\kappa(\beta)^{p-1}\lambda\bigr), 
\end{align*}
where for $c\in k$ and $\lambda=[\lambda_1:\lambda_2]\in\mathbb{P}^1(k)$ 
we set $c\lambda=[c\lambda_1:\lambda_2]$.
\end{prop}
\begin{proof}
Since 
$\mathcal{Z}\otimes\mathcal{X}_1^-\otimes\mathcal{X}_1^-\cong
\mathcal{Z}\otimes\mathcal{X}_1^+\cong \mathcal{Z}$ and 
$\mathcal{X}_1^-\otimes\mathcal{X}_1^-\otimes\mathcal{Z}\cong
\mathcal{X}_1^+\otimes\mathcal{Z}\cong \mathcal{Z}$
for 
any $\overline{U}$-module $\mathcal{Z}$, it is enough to show
\begin{align*}
\mathcal{E}_s^{+}(1;\lambda)\otimes\mathcal{X}_1^-&\cong
\mathcal{E}_s^{-}(1;-\lambda), \\
\mathcal{X}_1^-\otimes\mathcal{E}_s^{+}(1;\lambda)&\cong
\mathcal{E}_s^{-}\bigl(1;(-1)^{p-1}\lambda\bigr). 
\end{align*}
%

By Proposition \ref{prop:basis}, we can assume 
$\mathcal{X}_1^-=ka_0$, 
$\mathcal{E}_s^+(1;\lambda)
=\bigoplus_{n=0}^{s-1}kb_n\oplus\bigoplus_{m=0}^{p-s-1}kx_m$ with
$\overline{U}$-action given as that proposition.
Then $\mathcal{E}_s^{+}(1;\lambda)\otimes\mathcal{X}_1^-$ has 
basis 
$\{b_n\otimes a_0\}_{n=0,\ldots,s-1}\amalg
\{x_m\otimes a_0\}_{m=0,\ldots,p-s-1}$ and $\overline{U}$-action on 
these vectors is as follows:
\begin{align*}
K(b_n\otimes a_0)&=-q^{s-1-2n}b_n\otimes a_0, \quad
K(x_m\otimes a_0)=q^{p-s-1-2m}x_m\otimes a_0, \\
E(b_n\otimes a_0)&=\begin{cases}
-[n][s-n]b_{n-1}\otimes a_0 & (n\ne0)\\
-\lambda_2x_{p-s-1}\otimes a_0 & (n=0)
\end{cases}, \quad
E(x_m\otimes a_0)=\begin{cases}
[m][p-s-m]x_{m-1}\otimes a_0&(m\neq0)\\0&(m=0)
\end{cases},\\
F(b_n\otimes a_0)&=\begin{cases}
b_{n+1}\otimes a_0&(n\neq s-1)\\
\lambda_1x_0\otimes a_0&(n=s-1)\end{cases}, \quad
F(x_m\otimes a_0)=\begin{cases}
x_{m+1}\otimes a_0&(m\neq p-s-1)\\
0&(m=p-s-1)\end{cases}.
\end{align*}
This shows immediately 
$\mathcal{E}_s^{+}(1;\lambda)\otimes\mathcal{X}_1^-\cong
\mathcal{E}_s^{-}(1;-\lambda)$. 

Let us candider the second case. The module 
$\mathcal{X}_1^-\otimes\mathcal{E}_s^{+}(1;\lambda)$ has 
basis 
$\{(-1)^na_0\otimes b_n\}_{n=0,\ldots,s-1}\amalg
\{(-1)^ma_0\otimes x_m\}_{m=0,\ldots,p-s-1}$. In the following, we give 
explicit formulas of $\overline{U}$-action on 
these vectors. For simplicity, we denote
$\widetilde{b}_n=(-1)^na_0\otimes b_n$ and 
$\widetilde{x}_m=(-1)^ma_0\otimes x_m$.
\begin{align*}
K(\widetilde{b}_n)&=-q^{s-1-2n}\widetilde{b}_n, \quad
K(\widetilde{x}_m)=q^{p-s-1-2m}\widetilde{x}_m, \\
E(\widetilde{b}_n)&=\begin{cases}
-[n][s-n]\widetilde{b}_{n-1} & (n\ne0)\\
(-1)^{p-s-1}\lambda_2\widetilde{x}_{p-s-1} & (n=0)
\end{cases}, \quad
E(\widetilde{x}_m)=\begin{cases}
[m][p-s-m]\widetilde{x}_{m-1}&(m\neq0)\\0&(m=0)
\end{cases},\\
F(\widetilde{b}_n)&=\begin{cases}
\widetilde{b}_{n+1}&(n\neq s-1)\\
(-1)^s\lambda_1\widetilde{x}_0&(n=s-1)\end{cases}, \quad
F(\widetilde{x}_m)=\begin{cases}
\widetilde{x}_{m+1}&(m\neq p-s-1)\\
0&(m=p-s-1)\end{cases}.
\end{align*}
These formulas tells us 
$\mathcal{X}_1^-\otimes\mathcal{E}_s^{+}(1;\lambda)\cong
\mathcal{E}_s^{-}(1;\mu)$ with
$\mu=\bigl[(-1)^s\lambda_1:(-1)^{p-s-1}\lambda_2\bigr]=(-1)^{p-1}\lambda$. 
\end{proof}
 
Secondly let us compute 
$\mathcal{E}_s^{+}(1;\lambda)\otimes\mathcal{X}_2^+$ and
$\mathcal{X}_2^+\otimes\mathcal{E}_s^{+}(1;\lambda)$. 

\begin{lem}\label{lem:E generated}
Let $\mathcal{Z}$ be a $\overline{U}$-module and $s=1,\ldots,p$.

\noindent{\rm (i)}\ 
If $v\in \mathcal{Z}$ satisfies 
\[
 Kv=\pm q^{s-1}v, \quad 
 F^{p-1}v\neq0, \quad \text{and}\quad
 Ev=\alpha F^{p-1}v
\]
for some $s=1,\ldots,p-1$ and $\alpha\in k$, 
then $\bigoplus_{n=1}^{p-1}kF^nv$ is a submodule of $\mathcal{Z}$ 
isomorphic to $\mathcal{E}_s^{\pm}\bigl(1;[1:\alpha]\bigr)$.

\noindent{\rm (ii)}\ 
If $v\in \mathcal{Z}$ satisfies 
\[
 Kv=\pm q^{-s+1}v, \quad 
 E^{p-1}v\neq0, \quad \text{and}\quad
 Fv=0 
\]
for some $s=1,\ldots,p-1$, 
then $\bigoplus_{n=1}^{p-1}kE^nv$ is a submodule of $\mathcal{Z}$ 
isomorphic to $\mathcal{E}_s^{\pm}\bigl(1;[0:1]\bigr)$.

\noindent{\rm (iii)}\
If $s=p$ and $v\in \mathcal{Z}$ satisfies the conditions in 
{\rm (i)} or {\rm (ii)}, 
then $\bigoplus_{n=1}^{p-1}kF^nv$ or $\bigoplus_{n=1}^{p-1}kE^nv$, 
respectively, is a submodule of $\mathcal{Z}$ 
isomorphic to $\mathcal{X}_p^{\pm}$.
\end{lem}
\begin{proof}
The assertions follow by comparing the standard equations
\begin{align*}
 EF^n&=F^nE+[n]F^{n-1}\frac{q^{-n+1}K-q^{n-1}K^{-1}}{q-q^{-1}}, \\
 FE^n&=E^nF-[n]E^{n-1}\frac{q^{n-1}K-q^{-n+1}K^{-1}}{q-q^{-1}}
\end{align*} 
with Proposition \ref{prop:basis}. 
\end{proof}

\begin{prop}\label{prop:E and X_2^+}
For $s=1,\ldots,p-1$ and 
$\lambda=[\lambda_1:\lambda_2]\in\mathbb{P}^1(k)$ we have
\begin{align*}
\mathcal{E}_s^+(1;\lambda)\otimes\mathcal{X}_2^+&\cong
\mathcal{E}_{s-1}^+\biggl(1;\frac{[s]}{[s-1]}\lambda\biggr)\oplus
\mathcal{E}_{s+1}^+\biggl(1;\frac{[s]}{[s+1]}\lambda\biggr), \\
\mathcal{X}_2^+\otimes\mathcal{E}_s^+(1;\lambda)&\cong
\mathcal{E}_{s-1}^+\biggl(1;-\frac{[s]}{[s-1]}\lambda\biggr)\oplus
\mathcal{E}_{s+1}^+\biggl(1;-\frac{[s]}{[s+1]}\lambda\biggr), 
\end{align*}
where we put
$\mathcal{E}_{s-1}^+\bigl(1;\pm\frac{[s]}{[s-1]}\lambda\bigr)
=\mathcal{X}_p^-$ if $s=1$, and 
$\mathcal{E}_{s+1}^+\bigl(1;\pm\frac{[s]}{[s+1]}\lambda\bigr)
=\mathcal{X}_p^+$ if $s=p-1$. 
\end{prop}
\begin{proof}
It is enough to show that the modules on the left-hand sides 
have submodules isomorphic to direct summands on the right-hand sides, 
because any nonzero $\overline{U}$-module cannot be isomorphic to 
a submodule of 
$\mathcal{E}_{s-1}^+\bigl(1;\pm\frac{[s]}{[s-1]}\lambda\bigr)$ and 
$\mathcal{E}_{s+1}^+\bigl(1;\pm\frac{[s]}{[s+1]}\lambda\bigr)$
simultaneously.

As in the proof of Proposition \ref{prop:E and X_1^-}, 
we can take basis $\{b_n\otimes a_l\}\amalg
\{x_m\otimes a_l\}$ ($n=0,\ldots,s-1$, $m=0,\ldots,p-s-1$, $l=0,1$) of 
$\mathcal{E}_s^+(1;\lambda)\otimes\mathcal{X}_2^+$ in which
$\overline{U}$-action on $b_n$, $x_m$, $a_l$ is as 
Proposition \ref{prop:basis}. 
Let $v=[s]q^sb_0\otimes a_0+\lambda_2x_{p-s-1}\otimes a_1$. 
Then $Kv=q^sv$ and, using the standard equality
\[
 \Delta(F^n)=\sum_{k=0}^{n}q^{k(n-k)}\qbinom{n}{k}F^{n-k}K^{-k}\otimes F^{k}
\]
(where $\qbinom{n}{k}=\frac{[n]!}{[k]![n-k]!}$ with 
$[k]!=\prod_{l=1}^k [l]$), we have
\begin{align*}
 F^{p-1}v
&=[s]q^s\bigl(F^{p-1}b_0\otimes a_0
+q^{p-2}[p-1]F^{p-2}K^{-1}b_0\otimes Fa_0\bigr)\\
&=[s]\lambda_1(q^sx_{p-s-1}\otimes a_0-q^{-1}x_{p-s-2}\otimes a_1), \\
 Ev
&=[s]q^sEb_0\otimes Ka_0+\lambda_2(Ex_{p-s-1}\otimes Ka_1
+x_{p-s-1}\otimes Ea_1)\\
&=\bigl([s]q^{s+1}\lambda_2+\lambda_2\bigr)x_{p-s-1}\otimes a_0
-[p-s-1]q^{-1}\lambda_2x_{p-s-2}\otimes a_1\\
&=[s+1]\lambda_2(q^sx_{p-s-1}\otimes a_0-q^{-1}x_{p-s-2}\otimes a_1).
\end{align*}
Hence if $\lambda\ne[0:1]$, $v$ satisfies the condition 
of (i) (or (iii) when $s=p-1$) of the previous lemma. Therefore
$\mathcal{E}_s^+(1;\lambda)\otimes\mathcal{X}_2^+$ has a submodule 
isomorphic to 
$\mathcal{E}_{s+1}^+\bigl(1;\frac{[s]}{[s+1]}\lambda\bigr)$. 
If $\lambda=[0:1]$, 
one can verify that 
$v=q^sx_{p-s-1}\otimes a_0-q^{-1}x_{p-s-2}\otimes a_1$ satisfies
the condition of (ii) (or (iii) when $s=p-1$) 
of the previous lemma by using the equality 
\[
\Delta(E^n)=\sum_{k=0}^{n}q^{k(n-k)}\qbinom{n}{k}E^k\otimes E^{n-k}K^k. 
\]
In this case also
$\mathcal{E}_s^+(1;\lambda)\otimes\mathcal{X}_2^+$ has a submodule 
isomorphic to 
$\mathcal{E}_{s+1}^+\bigl(1;\frac{[s]}{[s+1]}\lambda\bigr)$. 
Similarly, let $w=b_1\otimes a_0-q[s-1]b_0\otimes a_1$, then we have
\begin{align*}
Kw=q^{s-2}w,\quad F^{p-1}w=-[s]\lambda_1x_{p-s-1}\otimes a_1, \quad
 Ew=-[s-1]\lambda_2x_{p-s-1}\otimes a_1.
\end{align*}
Hence the previous lemma shows that
$\mathcal{E}_s^+(1;\lambda)\otimes\mathcal{X}_2^+$ has a submodule 
isomorphic to 
$\mathcal{E}_{s-1}^+\bigl(1;\frac{[s]}{[s-1]}\lambda\bigr)$ for $\lambda\ne
[0:1]$. In the case of $\lambda=[0:1]$, let $w=x_{p-s-1}\otimes a_1$.
Then we have the same result by the similar argument. Consequently we have
$$\mathcal{E}_s^+(1;\lambda)\otimes\mathcal{X}_2^+\supset
\mathcal{E}_{s+1}^+\bigl(1;\frac{[s]}{[s+1]}\lambda\bigr)\oplus
\mathcal{E}_{s-1}^+\bigl(1;\frac{[s]}{[s-1]}\lambda\bigr).
$$ 
Since the dimension of each side is equal to $2p$, we have the first formula
of the proposition.\\ 

In the case of $\mathcal{X}_2^+\otimes\mathcal{E}_s^+(1;\lambda)$ 
one can prove the assertion by the same process:  
let $v=[s]a_0\otimes b_0+\lambda_2a_1\otimes x_{p-s-1}$ and 
$w=[s-1]a_1\otimes b_0-q^{s-1}a_0\otimes b_1$. Then
we have $Kv=q^sv$, $Kw=q^{s-2}w$ and
\[
F^{p-1}v=-[s]\lambda_1(qa_0\otimes x_{p-s-1}-a_1\otimes x_{p-s-2}),\quad
Ev=[s+1]\lambda_2(qa_0\otimes x_{p-s-1}-a_1\otimes x_{p-s-2}),
\]
\[
F^{p-1}v=-[s]\lambda_1a_1\otimes x_{p-s-1},\quad
Ev=[s-1]\lambda_2a_1\otimes x_{p-s-1}, 
\]
which leads us to the desired results.
\end{proof}

Thirdly, using Proposition \ref{prop:E and X_2^+}, 
we can calculate tensor products 
$\mathcal{E}_s^{+}(1;\lambda)\otimes\mathcal{X}_{s'}^{+}$ and
$\mathcal{X}_{s'}^{+}\otimes\mathcal{E}_s^{+}(1;\lambda)$ inductively 
on $s'$ as follows:
if $\mathcal{E}_s^{+}(1;\lambda)\otimes\mathcal{X}_{t}^{+}$ has 
known for $t\le s'-1$, the isomorphism
\begin{align*}
\bigl(\mathcal{E}_s^{+}(1;\lambda)\otimes\mathcal{X}_{s'-1}^{+}\bigr)
\otimes \mathcal{X}_2^+
&\cong 
\mathcal{E}_s^{+}(1;\lambda)\otimes(\mathcal{X}_{s'-1}^{+}
\otimes \mathcal{X}_2^+)
\cong
\mathcal{E}_s^{+}(1;\lambda)\otimes(\mathcal{X}_{s'-2}^{+}
\oplus\mathcal{X}_{s'}^+)\\
&\cong
\bigl(\mathcal{E}_s^{+}(1;\lambda)\otimes\mathcal{X}_{s'-2}^{+}\bigr)
\oplus
\bigl(\mathcal{E}_s^{+}(1;\lambda)\otimes\mathcal{X}_{s'}^+\bigr)
\end{align*}
determines the indecomposable decomposition of 
$\mathcal{E}_s^{+}(1;\lambda)\otimes\mathcal{X}_{s'}^+$. The explicit formulas 
are as follows:

\begin{prop}\label{prop:E and X_s^+}
For $s,s'=1,\ldots,p-1$ and $\lambda\in\mathbb{P}^1(k)$ 
we have 
\begin{align*}
\mathcal{E}_s^{+}(1;\lambda)\otimes\mathcal{X}_{s'}^{+}&\cong
\bigoplus_{t\in I_{s,s'}}\mathcal{E}_t^+\biggl(1;\frac{[s]}{[t]}\lambda\biggr)
\oplus \bigoplus_{t\in J_{s+s'}}\mathcal{P}_t^+
\oplus \bigoplus_{t\in J_{p-s+s'}}\mathcal{P}_t^-, \\
\mathcal{X}_{s'}^{+}\otimes\mathcal{E}_s^{+}(1;\lambda)&\cong
\bigoplus_{t\in I_{s,s'}}\mathcal{E}_t^+
\biggl(1;(-1)^{s'-1}\frac{[s]}{[t]}\lambda\biggr)
\oplus \bigoplus_{t\in J_{s+s'}}\mathcal{P}_t^+
\oplus \bigoplus_{t\in J_{p-s+s'}}\mathcal{P}_t^-.
\end{align*}
\end{prop}
\begin{proof}
Let us prove the first formula. 
From the exact sequence
\[
0\longrightarrow 
\mathcal{X}_{p-s}^{-}\otimes\mathcal{X}_{s'}^{+}
\longrightarrow
\mathcal{E}_s^{+}(1;\lambda)\otimes\mathcal{X}_{s'}^{+}
\longrightarrow
\mathcal{X}_s^{+}\otimes\mathcal{X}_{s'}^{+}
\longrightarrow0
\]
and the next decomposition formulas coming from 
Theorem \ref{thm:simple and simple} and 
Proposition \ref{prop:I and J}:
$$\mathcal{X}_{p-s}^{-}\otimes\mathcal{X}_{s'}^{+}
\cong \bigoplus_{t\in I_{s,s'}}\mathcal{X}_{p-t}^-
\oplus\bigoplus_{t\in J_{p-s+s'}}\mathcal{P}_t^-,\qquad
\mathcal{X}_{s}^{+}\otimes\mathcal{X}_{s'}^{+}
\cong \bigoplus_{t\in I_{s,s'}}\mathcal{X}_{t}^+
\oplus\bigoplus_{t\in J_{s+s'}}\mathcal{P}_t^+,$$ 
we have 
\[
\mathcal{E}_s^{+}(1;\lambda)\otimes\mathcal{X}_{s'}^{+}
\cong \bigoplus_{t\in I_{s,s'}}\mathcal{Z}_t
\oplus \bigoplus_{t\in J_{s+s'}}\mathcal{P}_t^+
\oplus \bigoplus_{t\in J_{p-s+s'}}\mathcal{P}_t^-,
\]
where $\mathcal{Z}_t$ is a nonprojective indecomposable module 
with an exact sequence
$0\longrightarrow \mathcal{X}_{p-t}^-
\longrightarrow \mathcal{Z}_t
\longrightarrow\mathcal{X}_t^+\longrightarrow0$ for 
each $t\in I_{s,s'}$. 

On the other hand, Proposition \ref{prop:E and X_2^+} and 
the calculation shown before the proposition, 
we see that a nonprojective indecomposable summand of   
$\mathcal{E}_s^{+}(1;\lambda)\otimes\mathcal{X}_{s'}^{+}$ 
must be of the form  
$\mathcal{E}_t^{+}\bigl(1;\frac{[s]}{[t]}\lambda\bigr)$ with
$t=1,\ldots,p-1$. 
Then we have 
$\mathcal{Z}_t\cong 
\mathcal{E}_t^{+}\bigl(1;\frac{[s]}{[t]}\lambda\bigr)$ since 
$\mathcal{Z}_t$ cannot be projective.  
Thus we have the first formula. 

The proof of the second formula is similar. 
\end{proof}

Finally, let us consider arbitrary cases. However, the result is an easy
consequence of Proposition \ref{prop:E and X_1^-} and \ref{prop:E and X_s^+}.
\begin{cor}
Let $\alpha,\beta\in\{+,-\}$. For $s,s'=1,\ldots,p-1$ and 
$\lambda\in\mathbb{P}^1(k)$ we have 
\begin{align*}
\mathcal{E}_s^{\alpha}(1;\lambda)\otimes\mathcal{X}_{s'}^{\beta}&\cong
\bigoplus_{t\in I_{s,s'}}\mathcal{E}_t^{\alpha\cdot\beta}
\biggl(1;\kappa(\beta)\frac{[s]}{[t]}\lambda\biggr)
\oplus \bigoplus_{t\in J_{s+s'}}\mathcal{P}_t^{\alpha\cdot\beta}
\oplus \bigoplus_{t\in J_{p-s+s'}}\mathcal{P}_t^{-\alpha\cdot\beta}, \\
\mathcal{X}_{s'}^{\beta}\otimes\mathcal{E}_s^{\alpha}(1;\lambda)&\cong
\bigoplus_{t\in I_{s,s'}}\mathcal{E}_t^{\alpha\cdot\beta}
\biggl(1;\kappa(\beta)^{s'-1}\frac{[s]}{[t]}\lambda\biggr)
\oplus \bigoplus_{t\in J_{s+s'}}\mathcal{P}_t^{\alpha\cdot\beta}
\oplus \bigoplus_{t\in J_{p-s+s'}}\mathcal{P}_t^{-\alpha\cdot\beta}.
\end{align*}
\end{cor}
\subsection{Rigidity}
For computing the remaining tensor products of indecomposable modules, 
we use a fact on finite-dimensional Hopf algebras. 

Let $A$ be a finite-dimensional Hopf algebra over $k$. Then it is known that
$A$-${\boldsymbol{\mathrm{mod}}}$ is a rigid tensor category 
({\it cf}. Appendix A). 

\begin{defn}
Let $\mathcal{Z}$ be a $A$-module. 
We define an $A$-module structure on 
the standard dual $D(\mathcal{Z})=\Hom_k(\mathcal{Z},k)$ by
$(a\varphi)(v)=\varphi\bigl(S(a)v\bigr)$ for
$a\in A$, $\varphi\in D(\mathcal{Z})$ and $v\in\mathcal{Z}$. 
\end{defn}

As a consequence of the rigidity, we have the following proposition which is
a central tool for computing tensor products ({\it cf}. Appendix A).
\begin{prop}\label{prop:rigidity}
For $A$-modules $\mathcal{Z}_1$, $\mathcal{Z}_2$, $\mathcal{Z}_3$ 
and $n\ge0$ we have
\[
 \Ext_A^n (\mathcal{Z}_1\otimes \mathcal{Z}_2,\mathcal{Z}_3)\cong
\Ext_A^n\bigl(\mathcal{Z}_1,\mathcal{Z}_3\otimes D(\mathcal{Z}_2)\bigr),\quad
 \Ext_A^n (\mathcal{Z}_1,\mathcal{Z}_2\otimes \mathcal{Z}_3)
\cong\Ext_A^n \bigl(D(\mathcal{Z}_2)\otimes \mathcal{Z}_1,\mathcal{Z}_3\bigr).
\]
\end{prop}

Let us compute $D(-)$ for our case $A=\overline{U}$.  

\begin{prop}\label{prop:D(X)}
For $s=1,\ldots,p-1$ and $\lambda\in\mathbb{P}^1(k)$ we have
\[
D(\mathcal{X}_s^{\pm})\cong\mathcal{X}_s^{\pm}, \quad
D\bigl(\mathcal{E}_s^+(1;\lambda)\bigr)\cong
\mathcal{E}_{p-s}^{-}\bigl(1;(-1)^s\lambda\bigr), \quad
D\bigl(\mathcal{E}_s^{-}(1;\lambda)\bigr)\cong
\mathcal{E}_{p-s}^{+}\bigl(1;(-1)^{p-s}\lambda\bigr). 
\]
\end{prop}
\begin{proof}
We only prove for $\mathcal{E}_s^{+}(1;\lambda)$. 
The other parts are similar. 

Take basis 
$\{b_n\}_{n=0,\ldots,s-1}\amalg\{x_m\}_{m=0,\ldots,p-s-1}$ of 
$\mathcal{E}_s^{+}(1;\lambda)$ as Proposition \ref{prop:basis}. 
Let 
$\{b_n^*\}_{n=0,\ldots,s-1}\amalg\{x_m^*\}_{m=0,\ldots,p-s-1}\subset 
D\bigl(\mathcal{E}_s^{+}(1;\lambda)\bigr)$ be 
the corresponding dual basis. Assume $\lambda\ne[0:1]$ and set $v=x_{p-s-1}^*$. 
Then we have
$$Kv=-q^{p-s-1}v,\quad F^{p-1}v=q^{(s-1)(p-1)}\lambda_1a_0^*, \quad
Ev=\mp q^{-s+1}\lambda_2a_0^*.$$
By Lemma \ref{lem:E generated} we have 
$D\bigl(\mathcal{E}_s^{+}(1;\lambda)\bigr)$ has a submodule
which is isomorphic to $\mathcal{E}_{p-s}^{-}\bigl(1;\mu\bigr)$, 
where $\mu=[q^{(s-1)(p-1)}\lambda_1:-q^{-s+1}\lambda_2]
=(-1)^s\lambda$. Thus we have the statement because
these modules have the same dimension. 

In the case of $\lambda=[0:1]$, the same argument as 
Proposition \ref{prop:E and X_2^+} is necessary. But we omit it in details. 
\end{proof}

\begin{prop}
For $s=1,\ldots,p-1$, $n\ge1$ and $\lambda\in\mathbb{P}^1(k)$ we have
\[
D\bigl(\mathcal{E}_s^{+}(n;\lambda)\bigr)\cong
\mathcal{E}_{p-s}^{-}\bigl(n;(-1)^s\lambda\bigr), \quad
D\bigl(\mathcal{E}_s^{-}(n;\lambda)\bigr)\cong
\mathcal{E}_{p-s}^{+}\bigl(n;(-1)^{p-s}\lambda\bigr). 
\]
\end{prop}
\begin{proof}
We prove the first formula, for the second one is proved similarly. 
Since $D$ preserves direct sum and dimension over $k$, 
we know that 
$D\bigl(\mathcal{E}_s^{+}(n;\lambda)\bigr)$ is an 
indecomposable module of dimension $pn$, 
therefore this is of the form $\mathcal{E}_t^{\pm}(n;\mu)$ or 
is projective (the latter case could occur only if $n\le2$). 

On the other hand, by Proposition \ref{prop:rigidity} we have
\begin{align*}
&\dim_k\Ext_{\overline{U}}^1
\bigl(D\bigl(\mathcal{E}_s^{+}(n;\lambda)\bigr),\mathcal{X}_s^+\bigr)\\
&\qquad=
\dim_k\Ext_{\overline{U}}^1
\bigl(D\bigl(\mathcal{E}_s^{+}(n;\lambda)\bigr)\otimes\mathcal{X}_1^+
,\mathcal{X}_s^+\bigr)
=\dim_k\Ext_{\overline{U}}^1
\bigl(\mathcal{X}_1^+,
\mathcal{E}_s^{+}(n;\lambda)\otimes\mathcal{X}_s^+\bigr)\\
&\qquad=
\dim_k\Ext_{\overline{U}}^1
\bigl(\mathcal{X}_1^+,
\mathcal{E}_s^{+}(n;\lambda)\otimes D(\mathcal{X}_s^+)\bigr)
=\dim_k\Ext_{\overline{U}}^1
\bigl(\mathcal{X}_1^+\otimes \mathcal{X}_s^+,
\mathcal{E}_s^{+}(n;\lambda)\bigr)\\
&\qquad=\dim_k\Ext_{\overline{U}}^1
\bigl(\mathcal{X}_s^+,
\mathcal{E}_s^{+}(n;\lambda)\bigr)=n, \\
&\dim_k\Ext_{\overline{U}}^1
\bigl(D\bigl(\mathcal{E}_s^{+}(n;\lambda)\bigr),\mathcal{E}_s^+(1;\mu)\bigr)\\
&\qquad=
\dim_k\Ext_{\overline{U}}^1
\bigl(D\bigl(\mathcal{E}_s^{+}(n;\lambda)\bigr)\otimes\mathcal{X}_1^+
,\mathcal{E}_s^+(1;\mu)\bigr)
=\dim_k\Ext_{\overline{U}}^1
\bigl(\mathcal{X}_1^+,
\mathcal{E}_s^{+}(n;\lambda)\otimes\mathcal{E}_s^+(1;\mu)\bigr)\\
&\qquad=\dim_k\Ext_{\overline{U}}^1
\Bigl(\mathcal{X}_1^+,
\mathcal{E}_s^{+}(n;\lambda)
\otimes D\bigl(\mathcal{E}_{p-s}^-\bigl(1;(-1)^{s}\mu\bigr)\bigr)\Bigr)\\
&\qquad=\dim_k\Ext_{\overline{U}}^1
\bigl(\mathcal{X}_1^+\otimes\mathcal{E}_{p-s}^-\bigl(1;(-1)^{s}\mu\bigr),
\mathcal{E}_s^{+}(n;\lambda)\bigr)
=\dim_k\Ext_{\overline{U}}^1
\bigl(\mathcal{E}_{p-s}^-\bigl(1;(-1)^{s}\mu\bigr),
\mathcal{E}_s^{+}(n;\lambda)\bigr)\\
&\qquad=\begin{cases}
1&((-1)^{s}\mu=-\lambda)\\
0&((-1)^{s}\mu\ne-\lambda)
\end{cases}.  
\end{align*}
Comparing these equalities with Proposition \ref{prop:ext} 
we have $D\bigl(\mathcal{E}_s^{+}(n;\lambda)\bigr)\cong
\mathcal{E}_{p-s}^{-}\bigl(n;(-1)^s\lambda\bigr)$ as desired. 
\end{proof}
\subsection{Tensor products of %
$\boldsymbol{\mathcal{E}_s^{\pm}(n;\lambda)}$ with %
simple modules}
Now we can calculate 
$\mathcal{E}_s^{\alpha}(n;\lambda)\otimes\mathcal{X}_{s'}^{\beta}$ and
$\mathcal{X}_{s'}^{\beta}\otimes\mathcal{E}_s^{\alpha}(n;\lambda)$ for 
general $n$ and $\alpha,\beta\in\{+,-\}$. However, by the similar method in 
Subsection 3.4, it is enough to consider the following cases; (a) 
$\alpha=\beta=+$ with arbitrary $s$ and $s'$, (b) $\beta=-$ and $s'=1$ with
arbitrary $\alpha$ and $s$. 

\begin{thm}\label{thm:E and X}
For $s,s'=1,\ldots,p-1$, $n\ge1$ and $\lambda\in\mathbb{P}^1(k)$ 
we have 
\begin{align*}
\mathcal{E}_s^{+}(n;\lambda)\otimes\mathcal{X}_{s'}^{+}&\cong
\bigoplus_{t\in I_{s,s'}}\mathcal{E}_t^+\biggl(n;\frac{[s]}{[t]}\lambda\biggr)
\oplus \bigoplus_{t\in J_{s+s'}}(\mathcal{P}_t^+)^n
\oplus \bigoplus_{t\in J_{p-s+s'}}(\mathcal{P}_t^-)^n, \\
\mathcal{X}_{s'}^{+}\otimes\mathcal{E}_s^{+}(n;\lambda)&\cong
\bigoplus_{t\in I_{s,s'}}\mathcal{E}_t^+
\biggl(n;(-1)^{s'-1}\frac{[s]}{[t]}\lambda\biggr)
\oplus \bigoplus_{t\in J_{s+s'}}(\mathcal{P}_t^+)^n
\oplus \bigoplus_{t\in J_{p-s+s'}}(\mathcal{P}_t^-)^n
\end{align*}
and 
\begin{align*}
\mathcal{E}_s^{\pm}(n;\lambda)\otimes\mathcal{X}_1^-&\cong
\mathcal{E}_s^{\mp}(n;-\lambda), \\
\mathcal{X}_1^-\otimes\mathcal{E}_s^{\pm}(n;\lambda)&\cong
\mathcal{E}_s^{\mp}\bigl(n;(-1)^{p-1}\lambda\bigr).  
\end{align*}
\end{thm}
\begin{proof}
We prove the first formula, for others are proved similarly. 
The same argument as Proposition \ref{prop:E and X_s^+} shows 
that there exists an isomorphism
\[
\mathcal{E}_s^{+}(n;\lambda)\otimes\mathcal{X}_{s'}^{+}
\cong \bigoplus_{t\in I_{s,s'}}\mathcal{Z}_t
\oplus \bigoplus_{t\in J_{s+s'}}(\mathcal{P}_t^+)^n
\oplus \bigoplus_{t\in J_{p-s+s'}}(\mathcal{P}_t^-)^n
\]
and an exact sequence
$0\longrightarrow (\mathcal{X}_{p-t}^-)^n
\longrightarrow \mathcal{Z}_t
\longrightarrow(\mathcal{X}_t^+)^n\longrightarrow0$ for 
each $t\in I_{s,s'}$. 
Moreover, by the exact sequence in Proposition 
\ref{prop:exact seq for E} and induction on $n$, we can assume that 
there exists an exact sequence
\[
0\longrightarrow 
\mathcal{E}_t^{+}\biggl(n-1;\frac{[s]}{[t]}\lambda\biggr)\longrightarrow 
\mathcal{Z}_t\longrightarrow 
\mathcal{E}_t^{+}\biggl(1;\frac{[s]}{[t]}\lambda\biggr)\longrightarrow 
0
\]
for each $t\in I_{s,s'}$. 

Let $t\in I_{s,s'}$. From
Proposition \ref{prop:rigidity}, Proposition \ref{prop:D(X)} 
and Proposition \ref{prop:ext} we have
\begin{align*}
&\dim_k\Ext_{\overline{U}}^1
\bigl(\mathcal{E}_s^{+}(n;\lambda)\otimes\mathcal{X}_{s'}^{+}
,\mathcal{X}_t^+\bigr)
\\&\qquad=
\dim_k\Ext_{\overline{U}}^1
\bigl(\mathcal{E}_s^{+}(n;\lambda)
,\mathcal{X}_t^+\otimes \mathcal{X}_{s'}^+\bigr)
=\dim_k\Ext_{\overline{U}}^1
\bigl(\mathcal{E}_s^{+}(n;\lambda)
,\mathcal{X}_s^+\bigr)=0, \\
&\dim_k\Ext_{\overline{U}}^1
\bigl(\mathcal{E}_s^{+}(n;\lambda)\otimes\mathcal{X}_{s'}^{+}
,\mathcal{X}_{p-t}^-\bigr)
\\&\qquad=
\dim_k\Ext_{\overline{U}}^1
\bigl(\mathcal{E}_s^{+}(n;\lambda)
,\mathcal{X}_{p-t}^-\otimes \mathcal{X}_{s'}^+\bigr)
=\dim_k\Ext_{\overline{U}}^1
\bigl(\mathcal{E}_s^{+}(n;\lambda)
,\mathcal{X}_{p-s}^-\bigr)=n, \\
&\dim_k\Ext_{\overline{U}}^1
\bigl(\mathcal{E}_s^{+}(n;\lambda)\otimes\mathcal{X}_{s'}^{+}
,\mathcal{E}_{t}^+(1;\mu)\bigr)
\\&\qquad=
\dim_k\Ext_{\overline{U}}^1
\bigl(\mathcal{E}_s^{+}(n;\lambda)
,\mathcal{E}_{t}^+(1;\mu)\otimes \mathcal{X}_{s'}^+\bigr)
=\dim_k\Ext_{\overline{U}}^1
\biggl(\mathcal{E}_s^{+}(n;\lambda)
,\mathcal{E}_{s}^+\biggl(1;\frac{[t]}{[s]}\mu\biggr)\biggr)\\
&\qquad=
\begin{cases}
1&\bigl(\lambda=\frac{[t]}{[s]}\mu\bigr)\\
0& \bigl(\lambda\ne\frac{[t]}{[s]}\mu\bigr)
\end{cases}.
\end{align*}
We note that 
$\mathcal{E}_s^{+}(n;\lambda)$ has no nontrivial first extension 
with modules from $\mathcal{C}(u)$ with $u\neq s$, 
and that $s\in I_{t,s'}$ by Proposition \ref{prop:I and J} (iv). 
This yields $\mathcal{Z}_t\cong 
\mathcal{E}_t^{+}\bigl(n;\frac{[s]}{[t]}\lambda\bigr)$ as desired. 
\end{proof}

Now we can calculate tensor products of 
$\mathcal{E}_s^{\alpha}(m;\lambda)$ with 
$\mathcal{M}_{s'}^{\beta}(n)$ or $\mathcal{W}_{s'}^{\beta}(n)$ by 
using projective covers and injective envelopes of 
$\mathcal{M}_{s'}^{\beta}(n)$, $\mathcal{W}_{s'}^{\beta}(n)$.
In the following, we only give the explicit formulas for $\alpha=\beta=+$, for 
simplicity. For the other combinations, we can easily calculate them by
the following theorem with the previous results.
The proof is analogous to that of Theorem \ref{thm:M and M}  
and is omitted. 

\begin{thm}\label{thm:E and M}
For $s,s'=1,\ldots,p-1$, $m\ge1$, $n\ge2$ and $\lambda\in\mathbb{P}^1(k)$ 
we have 
\begin{align*}
&\mathcal{E}_s^+(m;\lambda)\otimes\mathcal{M}_{s'}^+(n)\\
&\qquad\cong
\bigoplus_{t\in I_{s,s'}}
\biggl(\mathcal{E}_{p-t}^-\biggl(m;-\frac{[s]}{[t]}\lambda\biggr)
\oplus (\mathcal{P}_{t}^+)^{m(n-1)}\biggr)
\oplus\bigoplus_{t\in J_{s+s'}}(\mathcal{P}_t^+)^{m(n-1)}
\oplus\bigoplus_{t\in J_{2p-s-s'}}(\mathcal{P}_t^+)^{mn}\\
&\qquad\qquad\qquad\oplus
\bigoplus_{t\in J_{p+s-s'}}(\mathcal{P}_t^-)^{mn}\oplus
\bigoplus_{t\in J_{p-s+s'}}(\mathcal{P}_t^-)^{m(n-1)}, \\
&\mathcal{M}_{s'}^+(n)\otimes\mathcal{E}_s^+(m;\lambda)\\
&\qquad\cong
\bigoplus_{t\in I_{s,s'}}
\biggl(\mathcal{E}_{p-t}^-\biggl(m;(-1)^{s'}\frac{[s]}{[t]}\lambda\biggr)
\oplus (\mathcal{P}_{t}^+)^{m(n-1)}\biggr)
\oplus\bigoplus_{t\in J_{s+s'}}(\mathcal{P}_t^+)^{m(n-1)}
\oplus\bigoplus_{t\in J_{2p-s-s'}}(\mathcal{P}_t^+)^{mn}\\
&\qquad\qquad\qquad\oplus
\bigoplus_{t\in J_{p+s-s'}}(\mathcal{P}_t^-)^{mn}\oplus
\bigoplus_{t\in J_{p-s+s'}}(\mathcal{P}_t^-)^{m(n-1)}, \\
&\mathcal{E}_s^+(m;\lambda)\otimes\mathcal{W}_{s'}^+(n)\\
&\qquad\cong
\bigoplus_{t\in I_{s,s'}}
\biggl(\mathcal{E}_{t}^+\biggl(m;\frac{[s]}{[t]}\lambda\biggr)
\oplus (\mathcal{P}_{t}^+)^{m(n-1)}\biggr)
\oplus\bigoplus_{t\in J_{s+s'}}(\mathcal{P}_t^+)^{mn}
\oplus\bigoplus_{t\in J_{2p-s-s'}}(\mathcal{P}_t^+)^{m(n-1)}\\
&\qquad\qquad\qquad\oplus
\bigoplus_{t\in J_{p+s-s'}}(\mathcal{P}_t^-)^{m(n-1)}\oplus
\bigoplus_{t\in J_{p-s+s'}}(\mathcal{P}_t^-)^{mn}, \\
&\mathcal{W}_{s'}^+(n)\otimes\mathcal{E}_s^+(m;\lambda)\\
&\qquad\cong
\bigoplus_{t\in I_{s,s'}}
\biggl(\mathcal{E}_{t}^+\biggl(m;(-1)^{s'-1}\frac{[s]}{[t]}\lambda\biggr)
\oplus (\mathcal{P}_{t}^+)^{m(n-1)}\biggr)
\oplus\bigoplus_{t\in J_{s+s'}}(\mathcal{P}_t^+)^{mn}
\oplus\bigoplus_{t\in J_{2p-s-s'}}(\mathcal{P}_t^+)^{m(n-1)}\\
&\qquad\qquad\qquad\oplus
\bigoplus_{t\in J_{p+s-s'}}(\mathcal{P}_t^-)^{m(n-1)}\oplus
\bigoplus_{t\in J_{p-s+s'}}(\mathcal{P}_t^-)^{mn}. 
\end{align*}
\end{thm}
\subsection{Tensor products of %
$\boldsymbol{\mathcal{E}_s^{\pm}(m;\lambda)}$ and %
$\boldsymbol{\mathcal{E}_{s'}^{\pm}(n;\mu)}$}
As same as the second half of the previous subsection, we only calculate 
$\mathcal{E}_s^{+}(m;\lambda)
\otimes\mathcal{E}_{s'}^{+}(n;\mu)$. 

We note that there exist following exact sequences: 
\[
0\longrightarrow
\mathcal{E}_s^{+}(m;\lambda)\otimes(\mathcal{X}_{p-s'}^{-})^n
\longrightarrow  
\mathcal{E}_s^{+}(m;\lambda)\otimes\mathcal{E}_{s'}^{+}(n;\mu) 
\longrightarrow  
\mathcal{E}_s^{+}(m;\lambda)\otimes(\mathcal{X}_{s'}^{+})^n 
\longrightarrow0, 
\]
\[
0\longrightarrow
(\mathcal{X}_{p-s}^{-})^m\otimes\mathcal{E}_{s'}^{+}(n;\mu)
\longrightarrow  
\mathcal{E}_s^{+}(m;\lambda)\otimes\mathcal{E}_{s'}^{+}(n;\mu) 
\longrightarrow  
(\mathcal{X}_s^{+})^m\otimes\mathcal{E}_{s'}^{+}(n;\mu) 
\longrightarrow0.  
\]

The left and right terms of these sequences are 
computed by using Theorem \ref{thm:E and X}, 
which proves the next result: 

\begin{prop}\label{prop:E and E}
For $s,s'=1,\ldots,p-1$, $m,n\ge1$ and $\lambda,\mu\in\mathbb{P}^1(k)$ 
we have
\begin{align*}
&\mathcal{E}_s^+(m;\lambda)\otimes\mathcal{E}_{s'}^+(n;\mu)\\
&\qquad\cong
\bigoplus_{t\in I_{s,s'}}
\mathcal{V}_t(s,s';m,n;\lambda,\mu)
\oplus\bigoplus_{t\in J_{s+s'}}(\mathcal{P}_t^+)^{mn}
\oplus\bigoplus_{t\in J_{2p-s-s'}}(\mathcal{P}_t^+)^{mn}\\
&\qquad\qquad\qquad\oplus
\bigoplus_{t\in J_{p+s-s'}}(\mathcal{P}_t^-)^{mn}\oplus
\bigoplus_{t\in J_{p-s+s'}}(\mathcal{P}_t^-)^{mn}, 
\end{align*}
where $\mathcal{V}_t(s,s';m,n;\lambda,\mu)$ is a module 
in $\mathcal{C}(t)$. Moreover, there exist exact sequences
\[
0\longrightarrow
\mathcal{E}_{p-t}^-\biggl(m;-\frac{[s]}{[t]}\lambda\biggr)^n
\longrightarrow
\mathcal{V}_t(s,s';m,n;\lambda,\mu)
\longrightarrow
\mathcal{E}_t^+\biggl(m;\frac{[s]}{[t]}\lambda\biggr)^n
\longrightarrow 0, 
\]
\[
0\longrightarrow
\mathcal{E}_{p-t}^-\biggl(n;(-1)^s\frac{[s']}{[t]}\mu\biggr)^m
\longrightarrow
\mathcal{V}_t(s,s';m,n;\lambda,\mu)
\longrightarrow
\mathcal{E}_t^+\biggl(n;(-1)^{s-1}\frac{[s']}{[t]}\mu\biggr)^m
\longrightarrow 0. 
\]
\end{prop}

\bigskip

Let us determine the decomposition of 
$\mathcal{V}_t(s,s';m,n;\lambda,\mu)$ as a direct sum of indecomposable 
modules. 

\begin{thm}\label{thm:V_t}
For $s,s'=1,\ldots,p-1$, $t\in I_{s,s'}$, $m,n\ge1$, and
$\lambda,\mu\in\mathbb{P}^1(k)$ we have
\[
\mathcal{V}_t(s,s';m,n;\lambda,\mu)
\cong\begin{cases}
\mathcal{E}_t^+(l,\nu_t)
\oplus \mathcal{E}_{p-t}^-(l,-\nu_t) 
\oplus (\mathcal{P}_t^+)^{mn-l}
& \bigl(\frac{[s]}{[t]}\lambda=(-1)^{s-1}\frac{[s']}{[t]}\mu=\nu_t\bigr)\\[4pt]
(\mathcal{P}_t^+)^{mn}
& \bigl(\frac{[s]}{[t]}\lambda\ne(-1)^{s-1}\frac{[s']}{[t]}\mu\bigr)
\end{cases},  
\] 
where $l=\min\{m,n\}$. 
\end{thm}
\begin{proof}
We have
\begin{align*}
&\dim_k\Ext_{\overline{U}}^1
\bigl(\mathcal{E}_s^+(m;\lambda)\otimes\mathcal{E}_{s'}^+(n;\mu)
,\mathcal{X}_t^+\bigr)
\\&\qquad=\dim_k\Ext_{\overline{U}}^1
\bigl(\mathcal{E}_s^+(m;\lambda)
,\mathcal{X}_t^+\otimes\mathcal{E}_{p-s'}^-(n;(-1)^{s'}\mu)\bigr)\\
&\qquad=\dim_k\Ext_{\overline{U}}^1
\biggl(\mathcal{E}_s^+(m;\lambda)
,\mathcal{E}_{p-s}^-\biggl(n;(-1)^{s'+t-1}\frac{[s']}{[s]}\mu\biggr)\biggr)\\
&\qquad=
\begin{cases}
\min\{m,n\}&((-1)^{s-1}[s]\lambda=[s']\mu)\\
0& ((-1)^{s-1}[s]\lambda\ne[s']\mu)
\end{cases}, 
\qquad(\text{$t\equiv s-s'+1\ \bmod2$ for $t\in I_{s,s'}$})\\
&\dim_k\Ext_{\overline{U}}^1
\bigl(\mathcal{E}_s^+(m;\lambda)\otimes\mathcal{E}_{s'}^+(n;\mu)
,\mathcal{X}_{p-t}^-\bigr)
\\&\qquad=\dim_k\Ext_{\overline{U}}^1
\bigl(\mathcal{E}_s^+(m;\lambda)
,\mathcal{X}_{p-t}^-\otimes\mathcal{E}_{p-s'}^-(n;(-1)^{s'}\mu)\bigr)\\
&\qquad=\dim_k\Ext_{\overline{U}}^1
\biggl(\mathcal{E}_s^+(m;\lambda)
,\mathcal{E}_{s}^+\biggl(n;(-1)^{s'+t}\frac{[s']}{[s]}\mu\biggr)\biggr)\\
&\qquad=
\begin{cases}
\min\{m,n\}&((-1)^{s-1}[s]\lambda=[s']\mu)\\
0& ((-1)^{s-1}[s]\lambda\ne[s']\mu)
\end{cases}.  
\end{align*} 
These equalities show that, if $(-1)^{s-1}[s]\lambda\ne[s']\mu$,  
$\mathcal{V}_t(s,s';m,n;\lambda,\mu)$ is a projective module. 
Hence, by the exact sequences in the previous proposition,
it is isomorphic to $(\mathcal{P}_t^+)^{mn}$. 

From now on we assume $(-1)^{s-1}[s]\lambda=[s']\mu$. Set $\nu_t= 
\frac{[s]}{[t]}\lambda=(-1)^{s-1}\frac{[s']}{[t]}\mu$. Firstly assume $n=1$.
Then, from the equalities above, it is immediately to see that the 
nonprojective direct summand of $\mathcal{V}_t(s,s';m,1;\lambda,\mu)$ is 
isomorphic to $\mathcal{E}_t^+(1,\nu_t)\oplus \mathcal{E}_{p-t}^-(1,-\nu_t)$. 
Secondly, let us consider general cases. Using the result for $n=1$, we have
\begin{align*}
&\dim_k\Ext_{\overline{U}}^1
\bigl(\mathcal{E}_s^+(m;\lambda)\otimes\mathcal{E}_{s'}^+(n;\mu)
,\mathcal{E}_t^+(1;\nu_t)\bigr)
\\&\qquad=\dim_k\Ext_{\overline{U}}^1
\bigl(\mathcal{E}_s^+(m;\lambda)
,\mathcal{E}_t^+(1;\nu_t)\otimes\mathcal{E}_{p-s'}^-(n;(-1)^{s'}\mu)\bigr)
\\&\qquad=\dim_k\Ext_{\overline{U}}^1
\bigl(\mathcal{E}_s^+(m;\lambda)
,\mathcal{E}_s^+(1;\lambda)\oplus\mathcal{E}_{p-s}^-(1;-\lambda)\bigr)
\\&\qquad=2.  
\end{align*}
This equality and the previous equalities show that 
the nonprojective direct summand of 
$\mathcal{V}_t(s,s';m,n;\lambda,\mu)$ is 
isomorphic to 
$\mathcal{E}_t^+\bigl(\min\{m,n\},\nu_t\bigr)
\oplus \mathcal{E}_{p-t}^-\bigl(\min\{m,n\},-\nu_t\bigr)$. 
The assertion follows. 
\end{proof}

\vskip 5mm

Theorem \ref{thm:simple and simple},  
Theorem \ref{thm:projective and simple}, 
Corollary \ref{cor:projective and arbitrary}, 
Theorem \ref{thm:M and M}, 
Theorem \ref{thm:E and X}, 
Theorem \ref{thm:E and M}, 
Proposition \ref{prop:E and E}, 
Theorem \ref{thm:V_t} and obvious combination of them 
give indecomposable decomposition of 
tensor products of arbitrary $\overline{U}$-modules.\\

From the results in this section we have

\begin{prop}\label{prop;not-braided}
{\rm (i)}\ 
Let $\mathcal{Z}_1$, $\mathcal{Z}_2$ be 
$\overline{U}_q(\mathfrak{sl}_2)$-modules. 
If neither $\mathcal{Z}_1$ nor $\mathcal{Z}_2$ has any 
indecomposable summand of type $\mathcal{E}$, we have
$\mathcal{Z}_1\otimes\mathcal{Z}_2\cong\mathcal{Z}_2\otimes\mathcal{Z}_1$. 

\noindent{\rm (ii)}\ 
If $p=2$, for arbitrary $\overline{U}_q(\mathfrak{sl}_2)$-modules
$\mathcal{Z}_1$, $\mathcal{Z}_2$ we have
$\mathcal{Z}_1\otimes\mathcal{Z}_2\cong\mathcal{Z}_2\otimes\mathcal{Z}_1$. 

\noindent{\rm (iii)}\ 
If $p\ge3$, there exist $\overline{U}_q(\mathfrak{sl}_2)$-modules
$\mathcal{Z}_1$, $\mathcal{Z}_2$ such that 
$\mathcal{Z}_1\otimes\mathcal{Z}_2\not\cong
\mathcal{Z}_2\otimes\mathcal{Z}_1$. 
In particular, 
$\overline{U}_q(\mathfrak{sl}_2)$-${\boldsymbol{\mathrm{mod}}}$ is 
not a braided tensor category. 
\end{prop}
\begin{proof}
The assertions (i) and (ii) are clear. 
For (iii), set $\mathcal{Z}_1=\mathcal{E}_1^+\bigl(1;[1:1]\bigr)$ and 
$\mathcal{Z}_2=\mathcal{X}_2^+$. 
\end{proof}

As a by-product we have 
\begin{cor}
If $q$ is a primitive $2p$-th root of unity, $\overline{U}_q(\mathfrak{sl}_2)$ 
has no universal $R$-matrices for $p\geq 3$.
That is, it is not a quasi-triangular Hopf algebra.
\end{cor}

\begin{rem}
Let $\overline{U}_q^{\geq 0}$ be the $k$-subalgebra of 
$\overline{U}_q(\mathfrak{sl}_2)$ generated by $E,K,K^{-1}$. It is a 
$2p^2$-dimensional Hopf subalgebra of $\overline{U}_q(\mathfrak{sl}_2)$. 
By the quantum double construction, 
${\mathcal D}(\overline{U}_q^{\geq 0})\colon=D(\overline{U}_q^{\geq 0})\otimes
\overline{U}_q^{\geq 0}$ has a structure of a quasi-triangular Hopf algebra.
One can show that there is no surjective Hopf algebra homomorphism
${\mathcal D}(\overline{U}_q^{\geq 0})\to \overline{U}_q(\mathfrak{sl}_2)$.
This fact tells us $\overline{U}_q(\mathfrak{sl}_2)$ can not be obtained from
the usual quantum double construction, but it {\it does not} give a proof of
non-existence of universal $R$-matrices.
\end{rem}
\section{Complements}
%
%
%
\subsection{A quasi-triangular Hopf algebra $\overline{D}$}

The phenomenon which we showed in Proposition \ref{prop;not-braided} 
can be explained partly by considering a finite dimensional
Hopf $k$-algebra $\overline{D}$ which has a Hopf subalgebra 
isomorphic to $\overline{U}$. 
$\overline{D}$ is defined by generators $e$, $f$, $t$, $t^{-1}$ and 
relations 
\[
 tt^{-1}=t^{-1}t=1, \quad tet^{-1}=qe, \quad tft^{-1}=q^{-1}f, 
\]
\[
ef-fe=\frac{t^2-t^{-2}}{q-q^{-1}}, \quad
t^{4p}=1, \quad e^p=0, \quad f^p=0.
\]
The Hopf algebra structure on $\overline{D}$ is given by
\begin{align*}
\Delta&\colon e\longmapsto e\otimes t^2+1\otimes e, \quad
F\longmapsto f\otimes 1+t^{-2}\otimes f, \\
&\phantom{\colon}t\longmapsto t\otimes t, \quad 
t^{-1}\longmapsto t^{-1}\otimes t^{-1}, \\
\varepsilon&\colon e\longmapsto 0,\quad f\longmapsto 0, \quad
t\longmapsto 1, \quad t^{-1}\longmapsto 1, \\
S&\colon e\longmapsto -et^{-2}, \quad f\longmapsto -t^2f, \quad
t\longmapsto t^{-1}, \quad t^{-1}\longmapsto t. 
\end{align*}
$\overline{U}$ can be embedded into $\overline{D}$ as a Hopf subalgebra by 
\[\iota\colon
 E\longmapsto e, \quad F\longmapsto f, \quad K\longmapsto t^2. 
\]

We remark that finite-dimensional indecomposable $\overline{D}$-modules
are classified by Xiao (\cite{X2}, see also \cite{X1}, \cite{X3}). 
Those are parametrized by
the positive root system of type $A_3^{(1)}$ and some additional data.\\  

As in \cite{FGST1}, $\overline{D}$ is a quasi-triangular Hopf algebra
and has an universal $R$-matrix
\[
 \overline{R}=
\frac{1}{4p}\sum_{m=0}^{p-1}\sum_{n,j=0}^{4p-1}
\frac{(q-q^{-1})^m}{[m]!}q^{\frac{m(m-1)}{2}+m(n-j)-\frac{nj}{2}}
e^mt^n\otimes f^mt^j\in\overline{D}\otimes\overline{D}. 
\]
This shows that $\overline{D}$-${\boldsymbol{\mathrm{mod}}}$ is 
a braided tensor category. 

\begin{defn}
Let ${\mathcal Z}$ be a finite dimensional $\overline{U}$-module. The 
$\overline{U}$-action on $\mathcal Z$ is defined by a $k$-algebra 
homomorphism $\rho\colon\overline{U}\to \End_k({\mathcal Z})$. 
We call $\mathcal Z$ {\it liftable} if there
exists a $k$-algebra homomorphism $\rho'\colon\overline{D}\to 
\End_k({\mathcal Z})$ such that $\rho=\rho'\circ \iota$. The map 
$\rho'$ is called a {\it lifting} of $\rho$.
\end{defn}

The following lemma is easy to verify.
\begin{lem}
Each indecomposable $\overline{U}$-module except 
$\mathcal{E}_s^{\pm}(n;\lambda)$ $(\lambda\neq[1:0],[0:1])$ is liftable.
On the other hand, $\mathcal{E}_s^{\pm}(n;\lambda)$ $(\lambda\neq[1:0],[0:1])$
is not liftable. As a by-product, a universal $R$-matrix $\overline{R}$
can act on ${\mathcal Z}_1\otimes{\mathcal Z}_2$ for liftable modules
${\mathcal Z}_1$, ${\mathcal Z}_2$, and if either ${\mathcal Z}_1$ or
${\mathcal Z}_2$ is $\mathcal{E}_s^{\pm}(n;\lambda)$ $(\lambda\neq[1:0],[0:1])$,
$\overline{R}$ can not act on ${\mathcal Z}_1\otimes{\mathcal Z}_2$.
\end{lem}

As we already mentioned, Xiao \cite{X2} classify all finite-dimensional 
indecomposable $\overline{D}$-modules. In his list, there is the 
indecomposable $\overline{D}$-module $T^s(\alpha,\kappa,n)$ where 
$1\leq s\leq p-1$, $\alpha\in \{1,-1,\sqrt{-1},-\sqrt{-1}\}$, 
$\kappa=(\kappa_1,\kappa_2)\in (k^{\times})^2$ and $n$ is a positive integer. 
In Appendix B, we will give the explicit construction of 
$T^s(\alpha,\kappa,n)$. 

Assume $\alpha\in\{\pm 1\}$. As a $\overline{U}$-module, 
$T^s(\alpha,\kappa,n)$ decomposes into two indecomposable modules 
(for details, see Appendix B):
\[
T^s(\alpha,\kappa,n)\cong {\mathcal E}_s^+(n;\sqrt{\kappa_1\kappa_2}\,)
\oplus{\mathcal E}_s^+(n;-\sqrt{\kappa_1\kappa_2}\,).
\]
Here we set 
${\mathcal E}_s^+(n;\beta)\colon={\mathcal E}_s^+\bigl(n;[1:\beta]\bigr)$ for $\beta\in k$.

Let $\mathcal Z$ be a liftable $\overline{U}$-module and, by a fixed lifting
$\rho'\colon\overline{D}\to \End_k({\mathcal Z})$, we regard $\mathcal Z$ as
a $\overline{D}$-module. Since $\overline{D}$ has an universal $R$-matrix 
$\overline{R}$, there is an isomorphism of $\overline{D}$-modules:
\[
\sigma\overline{R}\colon T^s(\alpha,\kappa,n)\otimes{\mathcal Z}\overset{\sim}{\to}
{\mathcal Z}\otimes T^s(\alpha,\kappa,n),
\]
where $\sigma(a\otimes b)=b\otimes a$.
This isomorphism induces
\[
\bigl({\mathcal E}_s^+(n;\sqrt{\kappa_1\kappa_2}\,)\otimes{\mathcal Z}\bigr)
\oplus 
\bigl({\mathcal E}_s^+(n;-\sqrt{\kappa_1\kappa_2}\,)\otimes{\mathcal Z}\bigr)
\overset{\sim}{\to}
\bigl({\mathcal Z}\otimes{\mathcal E}_s^+(n;\sqrt{\kappa_1\kappa_2}\,)\bigr)
\oplus \bigl({\mathcal Z}\otimes
{\mathcal E}_s^+(n;-\sqrt{\kappa_1\kappa_2}\,)\bigr).
\]
Since $\overline{U}$ is a subalgebra of $\overline{D}$, the map above is also
an isomorphism of $\overline{U}$-modules. However, it interchanges the first 
and the second component, namely it induces isomorphisms of 
$\overline{U}$-modules
\[
{\mathcal E}_s^+(n;\sqrt{\kappa_1\kappa_2}\,)\otimes{\mathcal Z}
\overset{\sim}{\to} 
{\mathcal Z}\otimes{\mathcal E}_s^+(n;-\sqrt{\kappa_1\kappa_2}\,)\quad 
\text{and}\quad
{\mathcal E}_s^+(n;-\sqrt{\kappa_1\kappa_2}\,)\otimes{\mathcal Z}
\overset{\sim}{\to} 
{\mathcal Z}\otimes{\mathcal E}_s^+(n;\sqrt{\kappa_1\kappa_2}\,).
\]

This explains ``why'' 
the difference between $\mathcal{Z}_1\otimes\mathcal{Z}_2$ and 
$\mathcal{Z}_2\otimes\mathcal{Z}_1$ is no more than the sign differences 
in the parameters of the modules of type $\mathcal{E}^+$. For the case of type
$\mathcal{E}^-$, the situation is similar. 
\appendix
\section{Finite dimensional Hopf algebras}
In this appendix, we give a quick review on known results on representation
theory of finite dimensional Hopf algebras. These results can be found in
\cite{BK}, \cite{Ben}, \cite{CP}, \cite{K}, \cite{R}, and \cite{Sw}.
\subsection{Basic facts}
Let $\nk$ be a field and $A$ an algebra over $\nk$. For a right $A$-module $M$, 
the dual space $D(M)\colon=\Hom_{\nk}(M,\nk)$ has 
a left $A$-module structure defined by
\[
(a\rightharpoonup \lambda)(m)=\mu(ma)\quad (a\in A,\,\lambda\in D(M),
\,m\in M).
\]
Here we denote by $\rightharpoonup$ the left $A$-action on 
$D(M)$.\\

From now on we assume $A$ is a Hopf algebra with coproduct $\Delta$, 
counit $\varepsilon$ and antipode $S$. 
A {\it right integral} $\mu$ of $A$ is an element of $D(A)$ satisfying
\[
(\mu\otimes \id)\Delta(a)=\mu(a)1_A
\]
for all $a\in A$. Here $1_A$ is the unit of $A$. 
The following theorem is due to Sweedler \cite{Sw} (See also \cite{R}). 

\begin{thm}[\cite{Sw}]\label{thm:appendix}
Assume $A$ is a finite dimensional Hopf algebra over $\nk$. \\
{\rm (i)} Up to a scalar multiple, there uniquely exists a right integral 
$\mu$.\\
{\rm (ii)} Regarding $A$ as a right $A$-module, $D(A)$ has a 
left $A$-module structure. For a right integral $\mu$, the map $A\to D(A)$ 
defined by 
\[
a \mapsto (a\rightharpoonup\mu)
\]
is an isomorphism of left $A$-modules.\\
{\rm (iii)} $S$ is bijective.   
\end{thm}
\begin{rem}
The right integral of $\overline{U}_q({\mathfrak sl}_2)$ is given by
\[
\mu(F^iE^mK^n)=c\delta_{i,p-1}\delta_{m,p-1}\delta_{n,p+1}\quad 
(c\in k^{\times}). 
\]
\end{rem}
The following corollary follows from the second statement of the theorem.
\begin{cor}
If $A$ is a finite dimensional Hopf algebra, $A$ is a Frobenius algebra.
As a by-product, the following are equivalent:\\
{\rm (a)} $M$ is a projective $A$-module.\\
{\rm (b)} $M$ is an injective $A$-module. 
\end{cor} 
\subsection{Rigid tensor categories}
In this subsection, we introduce a notion of {\it rigid tensor categories}
following Bakalov and Kirillov, Jr. \cite{BK}. 

Let ${\mathcal C}$ be a monoidal category with the bifunctor 
$\otimes:{\mathcal C}\times{\mathcal C}\to{\mathcal C}$ and
the unit object ${\bf 1}\in \Ob{\mathcal C}$. For $V\in 
\Ob{\mathcal C}$, a {\it right dual} to $V$ is an object 
$D^R(V)$ with two morphisms
$$e_V^R\colon D^R(V)\otimes V \to {\bf 1},$$
$$i_V^R\colon{\bf 1}\to V\otimes D^R(V),$$
such that the two compositions
$$V\cong {\bf 1}\otimes V~\spmapright{i_V^R\otimes \id_V}~
V\otimes D^R(V)\otimes V~\spmapright{\id_V\otimes e_V^R}~
V\otimes {\bf 1}\cong V$$
and
$$D^R(V)\cong D^R(V)\otimes {\bf 1}
~\lspmapright{\id_{D^R(V)}\otimes i_V^R}~
D^R(V)\otimes V\otimes D^R(V)
~\lspmapright{e_V^R\otimes \id_{D^R(V)}}~
{\bf 1}\otimes D^R(V)\cong D^R(V)$$
are equal to $\id_V$ and $\id_{D^R(V)}$, respectively.

Similarly to the above, we define a {\it left dual} of $V$ to be an object
$D^L(V)$ with morphisms
$$e_V^L\colon V\otimes D^L(V) \to {\bf 1},$$
$$i_V^L\colon{\bf 1}\to D^L(V)\otimes V$$
and similar axioms.
\begin{defn}
A monoidal category ${\mathcal C}$ is called a {\it rigid tensor category} 
if every object in $\mathcal C$ has right and left duals.
\end{defn}
\begin{prop}\label{prop:hom}
Let $\mathcal C$ be a rigid tensor category and $V_1,V_2,V_3$ objects in 
$\mathcal C$.\\
{\rm (i)} $\Hom_{\mathcal C}(V_1,V_2\otimes V_3)\cong 
\Hom_{\mathcal C}(D^R(V_2)\otimes V_1,V_3)\cong 
\Hom_{\mathcal C}(V_1\otimes D^L(V_3),V_2)$.\\
{\rm (ii)} $\Hom_{\mathcal C}(V_1\otimes V_2,V_3)\cong
\Hom_{\mathcal C}(V_1,V_3\otimes D^R(V_2))\cong
\Hom_{\mathcal C}(V_2,D^L(V_1)\otimes V_3)$.
\end{prop}
\begin{proof} We only prove the first isomorphism of (ii). 
The others are proved by the similar way.

Define a map $\Phi\colon\Hom_{\mathcal C}(V_1\otimes V_2,V_3) \to 
\Hom_{\mathcal C}(V_1,V_3\otimes D^R(V_2))$ by
\vskip 1mm  
$$\Phi(f)\colon V_1\cong V_1\otimes \nk~\spmapright{\id\otimes i_{V_2}^R}~ 
V_1\otimes V_2\otimes D^R(V_2)
~\lspmapright{f\otimes \id_{D^R(V_2)}}~ V_3\otimes D^R(V_2)$$
for $f\in \Hom_{\mathcal C}(V_1\otimes V_2,V_3)$.  We remark that,
by the rigidity axioms, $\Phi(f)$ gives an element of 
$\Hom_{\mathcal C}(V_1,V_3\otimes D^R(V_2))$. Similarly we define
a well-defined map $\Psi\colon\Hom_{\mathcal C}(V_1,V_3\otimes D^R(V_2))\to 
\Hom_{\mathcal C}(V_1\otimes V_2,V_3)$ by
$$\Psi(g)\colon V_1\otimes V_2~\spmapright{g\otimes \id_{V_2}}~
V_3\otimes D^R(V_2)\otimes V_2~\spmapright{\id_{V_3}\otimes e_{V_2}^R}~
V_3\otimes\nk\cong V_3$$
for $g\in \Hom_{\mathcal C}(V_1,V_3\otimes D^R(V_2))$.

It is easy to see that $\Phi$ and $\Psi$ are inverse each other.
Thus, we have the statement.
\end{proof}
\subsection{The category of finite dimensional modules over a finite 
dimensional Hopf algebra}
Recall that $A$ is a finite dimensional Hopf algebra over a field $\nk$.
Let $A$-${\boldsymbol{\mathrm{mod}}}$ be the category of finite dimensional
left $A$-modules. It has a structure of a monoidal category associated with
the Hopf algebra structure of $A$. 

For a finite dimensional left $A$-module $V$, we define two left module
structure on $D(V)=\Hom_{\nk}(V,\nk)$: for $a\in A$, $\lambda\in D(V)$ 
and $v\in V$,
\begin{align*}
(a\rightharpoonup\lambda)(v)&=\lambda(S(a)v),\\ 
(a\rightharpoonup\lambda)(v)&=\lambda(S^{-1}(a)v).
\end{align*}
We denote by $D^R(V)$ the first left $A$-module
structure on $D(V)$ and by $D^L(V)$ the second one.
\begin{rem}
{(i)} Since the antipode $S$ is bijective 
(See Theorem \ref{thm:appendix} (iii)), $S^{-1}$ is a well-defined 
anti-isomorphism of $A$. However, $(A,\Delta,\varepsilon,S^{-1})$ 
is not a Hopf algebra in general.
More precisely $S^{-1}$ does not satisfy the axiom of an antipode.\\
{(ii)} If $S^2\ne\id_A$, $D^L(V)$ is not isomorphic to $D^R(V)$, 
in general. We remark that $S^2\ne\id_A$ for 
$A=\overline{U}_q(\mathfrak{sl}_2)$.
\end{rem}

By the construction, it is easy to see that
\[
D^R(D^L(V))\cong V\qquad \text{and}\qquad D^L(D^R(V))\cong V.
\]
The following proposition is easy to verify.

\begin{prop}\label{prop:rigid}
Let $V$ be an object in $A$-${\boldsymbol{\mathrm{mod}}}$, 
$\{v_i\}$ a basis of $V$ and $\{v_i^*\}$ the dual basis of $D(V)$.\\
{\em (i)} The $\nk$-linear maps 
$e_V^R\colon D^R(V)\otimes V\to\nk$ and 
$i_V^R\colon\nk\to V\otimes D^R(V)$ defined by
$$e_V^R(\lambda\otimes v)=\lambda(v)\qquad\text{and}\qquad
i_V^R(\alpha)=\alpha\left(\sum_i v_i\otimes v_i^*\right)$$
are homomorphisms of left $A$-modules, where we regard $\nk$ as a 
left $A$-module via the counit $\varepsilon$. 
Therefore $D^R(V)$ is the right dual to $V$.\\
{\em (ii)} Similarly, the $\nk$-linear maps 
$e_V^L\colon V\otimes D^L(V)\to \nk$ and 
$i_V^L\colon\nk\to D^L(V)\otimes V$ defined by
$$e_V^L(v\otimes\lambda)=\lambda(v)\qquad\text{and}\qquad
i_V^L(\alpha)=\alpha\left(\sum_i v_i^*\otimes v_i\right)$$
are homomorphisms of left $A$-modules. Therefore $D^L(V)$ is the
left dual to $V$.\\
{\em (iii)} $A$-${\boldsymbol{\mathrm{mod}}}$ is a rigid tensor category.
\end{prop}
As a consequence of the rigidity of $A$-${\boldsymbol{\mathrm{mod}}}$ and
Proposition \ref{prop:hom}, we have
\begin{cor}\label{cor:A-hom}
Let $V_1,V_2,V_3$ be objects in $A$-${\boldsymbol{\mathrm{mod}}}$.\\
{\em (i)} $\Hom_A(V_1,V_2\otimes V_3)\cong 
\Hom_A(D^R(V_2)\otimes V_1,V_3)\cong 
\Hom_A(V_1\otimes D^L(V_3),V_2)$.\\
{\em (ii)} $\Hom_A(V_1\otimes V_2,V_3)\cong
\Hom_A(V_1,V_3\otimes D^R(V_2))\cong
\Hom_A(V_2,D^L(V_1)\otimes V_3)$.
\end{cor}

\begin{cor}\label{cor:proj}
Let $P$ be a projective module. Then $P\otimes V$ and $V\otimes P$ are
also projective for any object $V$ in $A$-${\boldsymbol{\mathrm{mod}}}$.
\end{cor}
\begin{proof} We only show the projectivity of $P\otimes V$. Let $W_1$ and $W_2$
be objects in $A$-${\boldsymbol{\mathrm{mod}}}$, and $g\colon W_1\to W_2$ a
surjective $A$-homomorphism. It is enough to show that
$$g_*\colon\Hom_A(P\otimes V,W_1)\to \Hom_A(P\otimes V,W_2)$$
is surjective. 
Let us consider the following diagram:
$$\begin{array}{ccc}
\Hom_A(P\otimes V,W_1) &\spmapright{g_*}& \mbox{Hom}_A(P\otimes V,W_2)\\
\rmapdown{\wr} && \rmapdown{\wr}\\
\Hom_A(P,W_1\otimes D^R(V)) &\sbmapright{(g\otimes \id)_*}&
\Hom_A(P,W_2\otimes D^R(V))
\end{array}$$
\vskip 1mm
\noindent
where the vertical maps are the isomorphisms constructed in the proof of 
Proposition
\ref{prop:hom}. By the construction, this diagram is commutative. Since $P$ is
projective, we have $(g\otimes \id)_*$ is surjective. Thus, the map $g_*$ is 
also surjective.
\end{proof}

\begin{cor}[{\it cf}. Proposition \ref{prop:rigidity}]\label{cor:ext}
Let $V_1,V_2,V_3$ be objects in $A$-${\boldsymbol{\mathrm{mod}}}$. For any 
$n\geq 0$, we have the following.\\
{\em (i)} $\Ext^n_A(V_1,V_2\otimes V_3)\cong 
\Ext^n_A(D^R(V_2)\otimes V_1,V_3)\cong 
\Ext^n_A(V_1\otimes D^L(V_3),V_2)$.\\
{\em (ii)} $\Ext^n_A(V_1\otimes V_2,V_3)\cong
\Ext^n_A(V_1,V_3\otimes D^R(V_2))\cong
\Ext^n_A(V_2,D^L(V_1)\otimes V_3)$.
\end{cor}
\begin{proof}
We only prove the first isomorphism in (ii). 
Take a projective resolution of $V_1$:
$$\cdots \overset{d_2}{\to} P_1(V_1) \overset{d_1}{\to} P_0(V_1)
\overset{d_0}{\to} V_1 \to 0.$$
Then 
$$\Ext^n_A(V_1,V_3\otimes D^R(V_2))=
\frac{\Ker\bigl(d_{n+1}^*\colon\Hom_A(P_n(V_1),V_3\otimes D^R(V_2))
{\to} \Hom_A(P_{n+1}(V_1),V_3\otimes D^R(V_2))\bigr)}
{\Img\bigl(d_{n}^*\colon\Hom_A(P_{n-1}(V_1),V_3\otimes D^R(V_2))
{\to} \Hom_A(P_{n}(V_1),V_3\otimes D^R(V_2))\bigr)}.$$

Since $-\otimes V_2$ is an exact functor, the sequence
$$\cdots \spmapright{d_2\otimes\id_{V_2}}~ P_1(V_1)\otimes V_2 
\spmapright{d_1\otimes\id_{V_2}}~ P_0(V_1)\otimes V_2
\spmapright{d_0\otimes\id_{V_2}}~ V_1\otimes V_2 \to 0$$
is exact. Moreover, since $P_n(V_1)\otimes V_2$ is projective for any 
$n\geq 0$, this sequence gives a projective resolution of $V_1\otimes V_2$.
Therefore we have
$$
\Ext^n_A(V_1\otimes V_2,V_3)=
\frac{\Ker\bigl((d_{n+1}\otimes\id_{V_2})^*\colon
\Hom_A(P_n(V_1)\otimes V_2,V_3)\to
\Hom_A(P_{n+1}(V_1)\otimes V_2,V_3)\bigr)}
{\Img\bigl((d_{n}\otimes\id_{V_2})^*\colon
\Hom_A(P_{n-1}(V_1)\otimes V_2,V_3)\to
\Hom_A(P_{n}(V_1)\otimes V_2,V_3)\bigr)}.
$$
By the construction, there exists a commutative diagram:
$$\begin{array}{ccc}
\Hom_A(P_n(V_1),V_3\otimes D^R(V_2)) &\spmapright{d_{n+1}^*}& 
\Hom_A(P_{n+1}(V_1),V_3\otimes D^R(V_2))\\
\rmapdown{\wr} && \rmapdown{\wr}\\
\Hom_A(P_n(V_1)\otimes V_2,V_3) &
\sbmapright{(d_{n+1}\otimes \id)^*}&
\Hom_A(P_{n+1}(V_1)\otimes V_2,V_3)
\end{array}.$$
This diagram induces an isomorphism 
$\Ext^n_A(V_1,V_3\otimes D^R(V_2))\overset{\sim}{\to}
\Ext^n_A(V_1\otimes V_2,V_3)$.
\end{proof}
\section{The modules ${\mathcal E}^{+}_s(n;\lambda)$ and 
$T^s(\alpha,\kappa,n)$}
\subsection{The module ${\mathcal E}^{+}_s(n;\lambda)$}
Recall that ${\mathcal E}^{+}_s(n;\lambda)$ is defined as the image of
${\mathcal E}^{+}(n;\lambda)$ under the functor $\Phi_s$ where 
$1\leq s\leq p-1$ and $\lambda=[\lambda_1:\lambda_2]\in {\mathbb P}^1(k)$.
Since the explicit forms of primitive orthogonal idempotents of 
$\overline{U}$ are given by Arike \cite{Ari2}, one can determine the explicit 
structure of ${\mathcal E}^{+}_s(n;\lambda)$. 

The basis of ${\mathcal E}^{+}_s(n;\lambda)$ is 
$\bigl\{b_i^s(m),x_j^s(m)\bigm|0\leq i\leq s-1,0\leq j\leq p-s-1,
1\leq m\leq n\bigr\}$ and the action of $E,F,K^{\pm}$ are given as:
$$K^{\pm}b_i^s(m)=q^{\pm(s-1-2i)}b_i^s(m),\quad 
K^{\pm}x_j^s(m)=-q^{\pm(p-s-1-2j)}x_j^s(m),$$ 
$$Eb_i^s(m)=\begin{cases}
[i][s-i]b_{i-1}^s(m)& (i\ne 0),\\
\lambda_2x_{p-s-1}^s(m)+x_{p-s-1}^s(m-1)& (i=0),
\end{cases}
\quad 
Ex_j^s(m)=-[j][p-s-j]x_{j-1}^s(m),$$
$$Fb_i^s(m)=\begin{cases}
b_{i+1}^s(m)& (i\ne s-1),\\
\lambda_1x_{0}^s(m)&(i=s-1),
\end{cases}
\quad 
Fx_j^s(m)=x_{j+1}^s(m),$$
where we set $x_i^s(0)=0$ and $x_{p-s}^s(m)=0$.
\subsection{The module $T^s(\alpha,\kappa,n)$ and its decomposition as
$\overline{U}$-module}
Following Xiao \cite{X2}, let us
introduce the indecomposable $\overline{D}$-module $T^s(\alpha,\kappa,n)$ for
$1\leq s\leq p-1$, $\alpha\in 
\{1,-1,\sqrt{-1},-\sqrt{-1}\}$, $\kappa=(\kappa_1,\kappa_2)\in (k^{\times})^2$ 
and $n\in \mathbb{Z}_{>0}$. 
The basis of $T^s(\alpha,\kappa,n)$ is 
$\bigl\{{\bf e}_u^s(\alpha,m),\hat{\bf e}_u^s(\alpha,m)\bigm|
0\leq u\leq p-1,1\leq 
m\leq n\bigr\}$ and the action of $e,f,t^{\pm}$ is given as:
$$t^{\pm}{\bf e}_u^s(\alpha,m)=\alpha^{\pm}q^{\pm(s-1-2u)/2}{\bf e}_u^s(\alpha,m),
\quad
t^{\pm}\hat{\bf e}_u^s(\alpha,m)
=-\alpha^{\pm}q^{\pm(s-1-2u)/2}\hat{\bf e}_u^s(\alpha,m),
$$
$$e{\bf e}_u^s(\alpha,m)=
\begin{cases}
\alpha^2[u][s-u]{\bf e}_{u-1}^s(\alpha,m)& (u\ne 0),\\
\kappa_1\hat{\bf e}_{p-1}^s(\alpha,m)+\hat{\bf e}_{p-1}^s(\alpha,m-1)& (u=0),
\end{cases}$$
$$e\hat{\bf e}_u^s(\alpha,m)=
\begin{cases}
\alpha^2[u][s-u]\hat{\bf e}_{u-1}^s(\alpha,m)& (u\ne 0),\\
\kappa_2{\bf e}_{p-1}^s(\alpha,m)+{\bf e}_{p-1}^s(\alpha,m-1)& (u=0),
\end{cases}$$
$$f{\bf e}_u^s(\alpha,m)={\bf e}_{u+1}^s(\alpha,m),\quad
f\hat{\bf e}_u^s(\alpha,m)=\hat{\bf e}_{u+1}^s(\alpha,m),$$
where ${\bf e}_u^s(\alpha,0)=\hat{\bf e}_u^s(\alpha,0)=0$ and 
${\bf e}_p^s(\alpha,m)=\hat{\bf e}_p^s(\alpha,m)=0$.\\

Assume $\alpha^2=1$. Consider an invertible $(2n\times 2n)$ matrix $Q$ 
which satisfies
$$Q^{-1}\left(\begin{array}{cc}
O & J(n;\kappa_2)\\
J(n;\kappa_1) & O
\end{array}\right)Q=\left(\begin{array}{cc}
J(n;\sqrt{\kappa_1\kappa_2}\,) & O\\
O & J(n;-\sqrt{\kappa_1\kappa_2}\,)
\end{array}\right)$$
where $J(n;\beta)$ is the $(n\times n)$-Jordan cell with the eigenvalue 
$\beta$.
Define ${\bf b}_u^s(\alpha,m),\hat{\bf b}_u^s(\alpha,m)$ 
$(0\leq u\leq p-1,1\leq m\leq n)$ by
$$({\bf b}_u^s(\alpha,1),\cdots,{\bf b}_u^s(\alpha,n),\hat{\bf b}_u^s(\alpha,1),
\cdots,\hat{\bf b}_u^s(\alpha,n)):=
({\bf e}_u^s(\alpha,1),\cdots,{\bf e}_u^s(\alpha,n),\hat{\bf e}_u^s(\alpha,1),
\cdots,\hat{\bf e}_u^s(\alpha,n))Q$$
and a $k$-linear isomorphism $\Psi:T^s(\alpha,\kappa,n)\to 
{\mathcal E}_s^+(n;\sqrt{\kappa_1\kappa_2})\oplus
{\mathcal E}_s^+(n;-\sqrt{\kappa_1\kappa_2})$ by
$${\bf b}_u^s(\alpha,m)\mapsto
\begin{cases}
b_u^{s,+}(m)& (0\leq u\leq s-1),\\
x_{u-s}^{s,+}(m)& (s\leq u\leq p-1),
\end{cases}
\quad
\hat{\bf b}_u^s(\alpha,m)\mapsto
\begin{cases}
b_u^{s,-}(m)& (0\leq u\leq s-1),\\
x_{u-s}^{s,-}(m)& (s\leq u\leq p-1),
\end{cases}
$$
where we denote by $\{b_i^{s,\pm}(m),x_j^{s,\pm}(m)\}$
the basis of ${\mathcal E}_s^+(n;\pm\sqrt{\kappa_1\kappa_2}\,)$
which is introduced in the previous subsection. By the construction, it is 
easy to see that $\Psi$ is an isomorphism of $\overline{U}$-modules.
 

\end{document}